\documentclass{amsart}

\usepackage{amsthm, amsmath, amssymb, xypic, url, verbatim,graphics}
\usepackage{algorithmic}
\usepackage{algorithm}
\usepackage{lscape, multirow}
\usepackage[curve]{xy}
\usepackage{indent}

\numberwithin{equation}{section}
\numberwithin{table}{section}

\newtheorem{theorem}{Theorem}[section]

\theoremstyle{definition}
\newtheorem{definition}[theorem]{Definition}
\newtheorem{remark}[theorem]{Remark}
\newcommand{\C}{C}
\newcommand{\F}{\mathbb{F}}
\newcommand{\G}{\mathbb{G}}
\newcommand{\N}{\mathbb{N}}
\newcommand{\Q}{\mathbb{Q}}

\newcommand{\Z}{\mathbb{Z}}
\newcommand{\Jac}{\textrm{Jac}}
\newcommand{\Gal}{\textrm{Gal}}
\newcommand{\inv}{\textrm{Inv}}
\newcommand{\ra}{\rightarrow}
\newcommand{\norm}{\textrm{norm}}
\newcommand{\lc}{\textrm{lc}}
\renewcommand{\pmod}[1]{\ \left({\rm mod\ } #1 \right)}
\renewcommand{\th}{\textrm{th}}

\DeclareMathOperator{\Tr}{Tr}

\DeclareMathOperator{\jac}{Jac}

\DeclareMathOperator{\Aut}{Aut}

\newcommand{\Div}{\operatorname{Div}}
\newcommand{\Ppl}{\operatorname{Ppl}}
\newcommand{\Pic}{\operatorname{Pic}}
\newcommand{\ord}{\operatorname{ord}}
\renewcommand{\div}{\operatorname{div}}
\newcommand{\barK}{\overline{K}}

\newenvironment{stage}[1]{\noindent \textsc{Stage #1}:\begin{indentation}{1.5em}{0em}}{\end{indentation}}

\begin{document}
\bibliographystyle{amsplain}

%%%%%%%%%%%%%%%%%%%%%%%%%%
%%%%%%begin frontmatter
%%%%%%%%%%%%%%%%%%%%%%%%%%
\title{Pairings on Hyperelliptic Curves}
\author[J.\ Balakrishnan]{Jennifer Balakrishnan}
\address{Dept.\ of Mathematics\\
Massachusetts Institute of Technology\\
Cambridge, MA 02139}
\email{jen@math.mit.edu}
\author[J.\ Belding]{Juliana Belding}
\address{Dept.\ of Mathematics\\
Harvard University \\
Cambridge, MA 02138}
\email{jbelding@math.harvard.edu}
\author[S.\ Chisholm]{Sarah Chisholm}
\address{Dept.\ of Mathematics and Statistics\\
        University of Calgary\\
        Calgary, Alberta, Canada T2N 1N4}
\email{chisholm@math.ucalgary.ca}
\author[K.\ Eisentr\"ager]{Kirsten Eisentr\"ager}
\address{Dept.\ of Mathematics\\
The Pennsylvania State University\\
University Park, PA 16802}
\email{eisentra@math.psu.eu}
\author[K.\ Stange]{Katherine E. Stange}
\address{Dept.\ of Mathematics\\
Simon Fraser University\\
Burnaby, British Columbia, Canada V5A 1S6}
\email{kestange@sfu.ca}
\author[E.\ Teske]{Edlyn Teske}
\address{Dept.\ of Combinatorics and Optimization\\
        University of Waterloo\\
        Waterloo, Ontario, Canada N2L 3G1}
\email{eteske@uwaterloo.ca}

\dedicatory{Dedicated to the memory of Isabelle D\'ech\`ene (1974-2009)}

\date{\today}
\begin{abstract}
  We assemble and reorganize the recent work in the area of
  hyperelliptic pairings: We survey the research on constructing
  hyperelliptic curves suitable for pairing-based cryptography. We
  also showcase the hyperelliptic pairings proposed to date, and develop a
  unifying framework. We discuss the techniques used to optimize the
  pairing computation on hyperelliptic curves, and present many
  directions for further research.
\end{abstract}

\maketitle
%%%%%%%%%%%%%%%%%%%%%%%%%%%%%%%%%%%%
%%%%%%%%%%%%%%% end frontmatter
%%%%%%%%%%%%%%%%%%%%%%%%%%%%%%%%%%%%

%****************************************************************************

\section{Introduction}
\label{s:intro}

Numerous cryptographic protocols for secure key exchange and digital
signatures are based on the computational infeasibility of the
discrete logarithm problem in the underlying group.  Here, the most
common groups in use are multiplicative groups of finite fields and
groups of points on elliptic curves over finite fields.  Over the past
years, many new and exciting cryptographic schemes based on pairings
have been suggested, including one-round three-way key establishment,
identity-based encryption, and short signatures
\cite{BonehFranklin:2003,BLS:2004,Joux:2004,PatersonCh:2005}.
Originally, the Weil and Tate (-Lichtenbaum) pairings on supersingular elliptic
curves were proposed for such applications, providing non-degenerate
bilinear maps that are efficient to evaluate. Over time potentially more efficient
pairings have been found, such as the eta \cite{BGOS:2007}, Ate
\cite{HSV:2006} and R-ate \cite{LLP:2009} pairings.  Computing any
of these pairings involves finding functions with prescribed zeros and
poles on the curve, and evaluating those functions at divisors.

As an alternative to elliptic curve groups, Koblitz \cite{Koblitz:89} 
suggested Jacobians of hyperelliptic curves for use in cryptography.
In particular, hyperelliptic curves of low genus represent a competitive
choice. %\cite{lange/paar/pelzl}.
In 2007, Galbraith, Hess and Vercauteren \cite{GHV:2007} summarized
the research on hyperelliptic pairings to date and compared the
efficiency of pairing computations on elliptic and hyperelliptic curves. In this rapidly moving
area, there have been several new developments since their survey: First, new
pairings have been developed for the elliptic case, including
so-called optimal pairings by Vercauteren \cite{Vercauteren:2008} and a framework for
elliptic pairings by Hess \cite{Hess:2008}. Second, several constructions of
ordinary hyperelliptic curves suitable for pairing-based cryptography
have been found \cite{Freeman:2008, FSS:2008, Satoh:2009, FreemanSatoh}.

In this paper, we survey
\begin{itemize}
\item the constructions of hyperelliptic curves suitable for pairings, especially in the ordinary case, %(Section \ref{s:friendly})
\item the hyperelliptic pairings proposed to date, %(Section \ref{s:hecpairings}); 
and
\item the techniques to optimize computations of hyperelliptic pairings. %(Section \ref{s:fast}. 
\end{itemize} 
We also
\begin{itemize}
\item give a unifying framework for hyperelliptic pairings which includes many of the recent variations of the Ate pairing, % (Sections \ref{ss:framework}, \ref{ss:examples}) ; 
and
\item present a host of potential further improvements.% (Section \ref{s:faster}).
\end{itemize}
In this paper, we do {\em not} provide any comparative implementation, or give recommendations on which pairings should be used to satisfy certain user-determined criteria; this is left for future work. 

In our presentation, we focus on the case of genus $2$ hyperelliptic
curves and their Jacobians. Among all curves of higher genus, such
curves are of primary interest for cryptographic applications: On the
one hand, we find explicit formulae along with various optimizations
(e.g., \cite{Lange:2005, WPP:2005}), providing for an arithmetic that
is somewhat competitive with elliptic curves. On the other hand,
the security is exactly the same as in the elliptic case, with the
best attacks on the discrete logarithm problem in the Jacobian being
square-root attacks based on the Pollard rho method
(cf. \cite{FreyLangeCh23:2006}).
However,
Galbraith, Hess and Vercauteren \cite[\S10.1]{GHV:2007} argue that {\em pairing} computations on
hyperelliptic curves will always be slow compared to
elliptic curves: The most expensive part of a standard
Tate pairing computation consists of repeatedly evaluating some function on a
divisor and computing the product of the values obtained. Both in the elliptic
and in the hyperelliptic case these divisors are defined over fields of
the same size, but the functions in the hyperelliptic case are more
complicated.

\begin{figure}[t]
\caption{Classification of hyperelliptic pairings}
\label{tab:tree}
\begin{equation*}\xymatrix{ & &
\textbf{\large{\txt{Hyperelliptic\\pairings}}} \ar@{-}[dll] \ar@{-}[dl]
\ar@{-}[d] \ar@{-}[dr] &  \\
\txt{Weil}  & \txt{Tate} & \txt{Modified\\Tate} & \txt{Ate\\family}
\ar@{-}[dll] \ar@{-}[dl] \ar@{-}[d] \\
& \txt{Eta \\ (twisted \\only)} & \txt{Ate\\(no final \\exponent)} &
\txt{Hess-\\Vercauteren\\(HV)} \ar@{-}[dl] \ar@{-}[d] \\
&&  \txt{Ate$_i$}  & \txt{R-ate}}
\end{equation*}
\end{figure}
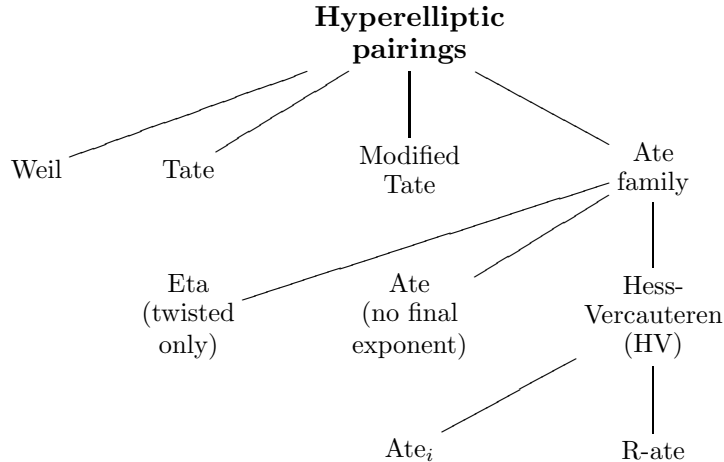

Figure \ref{tab:tree} represents the collection of hyperelliptic pairings at a glance. 
For use in pairing-based applications, originally the Weil and Tate pairings 
were proposed. 
The Weil pairing is much more expensive to
compute than the Tate pairing, so it is not used in practice. 
The
pairings in the Ate family are potentially more efficient than the
Tate pairing. 
Historically, the eta pairing was the first pairing to shorten the  
length of the Miller loop. It is defined on supersingular curves only 
and requires a final exponentiation. It gave rise to the Ate pairings  
which are defined for all curves. 
The hyperelliptic Ate pairing (which has a different definition than the elliptic Ate pairing\,!) has the advantage that its
loop length is roughly half of the length of the Miller loop for the
Tate pairing. It also is special in that it requires no final exponentiation (while the elliptic Ate pairing does require one). 
Other variations of the
Ate pairing include the Hess-Vercauteren (HV) pairings. These are 
the pairings captured by our unifying framework, which generalizes work for the elliptic case by Hess \cite{Hess:2008} and Vercauteren \cite{Vercauteren:2008}. 
HV pairings also have potentially shorter Miller loops than the Ate pairing, depending on the embedding degree  of the Jacobian. All of the HV pairings
involve a
final exponentiation. Two examples of HV pairings are the R-ate and
the Ate$_i$ pairings.
Table~\ref{pairingtable} in Section~\ref{s:fast} gives more details about
the differences and merits of each pairing.

Our paper is organized as follows.
In Section~\ref{s:math} we review some of the background on Jacobians
of hyperelliptic curves. Section~\ref{s:friendly} discusses
hyperelliptic curves of low embedding degree and what is known about
constructing them.  Section~\ref{s:pairings} gives an overview of the
different pairings on hyperelliptic curves 
following the classification in Figure \ref{tab:tree}.
We also introduce the HV pairing framework,  give a direct
proof of the non-degeneracy and bilinearity of the pairings captured by this framework and discuss how the Ate and R-ate pairings fit in.
Section~\ref{s:fast} describes the adaptation of Miller's algorithm to
the hyperelliptic setting, presents common optimizations and
compares all pairings according to their key characteristics of loop
length and final exponentiation.
Section~\ref{s:faster} presents numerous problems for future work.

%*******************************************************

\section{Jacobians of Hyperelliptic Curves}
\label{s:math}
In this section, we fix some notation and terminology that will be used throughout the paper.

\subsection{Hyperelliptic curves}\label{hyperellipticcurvedefinition}
A \emph{hyperelliptic curve} $\C$ over a field $K$ is a non-singular projective curve of the form
\[
\C: y^2 + H(x)y = F(x) \in K[x,y].
\]
Let $g$ be the genus of the curve.  Throughout this paper, we restrict to the case where $F$ is monic, $\deg F(x) = 2g+1$, and $\deg H(x) \leq g$, so that $\C$ has one point at infinity, denoted $P_{\infty}$. When $g = 1$, $\C$ is an \emph{elliptic curve}. For significant parts of our discussion, we will consider the case where $g = 2$.

Although the points of a genus $g \geq 2$ hyperelliptic curve do not form a group, there is an involution of the curve taking $P = (x,y)$ to the point $(x,-y-H(x)),$ which we will denote $-P$.  Then, in accordance with the notation, $-(-P)=P$.

\subsection{Divisors and abelian varieties}

Let $K$ be a field over which $\C$ is defined, and let $\barK$ its algebraic closure.  A \emph{divisor} $D$ on the curve $\C$ is a formal sum over all symbols $(P)$, where $P$ is a $\barK$-point of the curve:
\[
D=\sum_{P \in \C(\barK)} n_P(P),
\]
where all but finitely many of the coefficients $n_P \in \mathbb{Z}$ are zero.  The collection of divisors forms an abelian group $\Div(\C)$.  The \emph{degree} of a divisor is the sum
\[
\sum_{P \in \C(\barK)}n_P \in \mathbb{Z}, 
\]
and the {\em support} of a divisor is the set of points of the divisor
with non-zero coefficients $n_P$. For any rational function $f$ on $\C$, there is an associated divisor
\[
\div(f) = \sum_{P \in \C(\barK)} \ord_P(f)(P)
\]
which encodes the number and location of its zeroes and poles.  Any divisor which is the divisor of a function in this way is called a \emph{principal divisor}.

An element $\sigma$ in the Galois group of $\barK$ over $K,$ $\Gal(\barK/K),$ acts on a divisor as follows:
\[
\Big(\sum_{P \in \C(\barK)}n_P(P)\Big)^\sigma =
\sum_{P \in \C(\barK)}n_P(P^\sigma).
\]
In particular, let $L$ be any intermediate field $K \subset L \subset \barK.$ Consider a function $f$ defined over $L;$ then $\div(f)$ is fixed by elements of $\Gal(\barK/L)$.  In fact, $\div(f)^\sigma = \div(f^\sigma)$.

We give names to various sets of collections: $\Div(\C)$ of divisors, $\Div^0(\C)$ of degree zero divisors, $\Ppl(\C)$ of principal divisors, $\Div_K(\C)$ of divisors invariant under the action of $\Gal(\barK /K)$,  $\Div_K^0(\C)$ of degree zero divisors invariant under the action of $\Gal(\barK /K)$, and $\Ppl_K(\C)$ of principal divisors invariant under the action of $\Gal(\barK /K)$.

These are all abelian groups, which have the following subgroup relations:

\[
\begin{array}{ccccc}
\Div(C) & \supset  & \Div^0(C) & \supset  & \Ppl(C) \\[0.5em]
 \cup & & \cup & & \cup \\[0.5em]
\Div_K(C) & \supset & \Div^0_K(C)  & \supset & \Ppl_K(C). \\
\end{array}
\]

We make note of certain quotient groups:
\begin{align*}
\Pic(\C) &:= \Div(\C) / \Ppl(\C), & \Pic^0(\C) &:= \Div^0(\C) / \Ppl(\C), \\
\Pic_K(\C) &:= \Div_K(\C) / \Ppl_K(\C), & \Pic^0_K(\C) &:= \Div^0_K(\C) / \Ppl_K(\C).
\end{align*}

Elements of these quotient groups are equivalence classes of divisors. Divisors $D_1$ and $D_2$ of the same class are said to be \emph{linearly equivalent}, and we write $D_1 \sim D_2$.

Recall that an elliptic curve is an example of an abelian variety.  In general, an \emph{abelian variety} $A$ over $K$ is a projective algebraic variety over $K$ along with a group law $\varphi: A \times A \rightarrow A$ and an inverse map $\inv: A \rightarrow A$ sending $x \mapsto x^{-1}$ such that $\varphi$ and $\inv$ are morphisms of varieties, both defined over $K$.

For an abelian variety $A$, a field $K$ and an integer $r$, we let $A(K)[r]$ denote the set of $r$-torsion points of $A$ defined over $K$, that is, the set of points in $A(K)$ of order dividing $r$. Now suppose $A$ is an abelian variety over $\F_q$, with $q = p^m$. We say that $A$ is \emph{simple} if it is not
isogenous over $\F_q$ to a product of lower dimensional abelian varieties. We call $A$ {\em absolutely simple} if it is simple over $\overline{\F}_q$. We say $A$ is \emph{supersingular} if $A$ is isogenous over $\overline{\F}_q$ to a power of a supersingular elliptic curve. (An elliptic curve $E$ is supersingular if $E(\overline{\F}_q)$ has no
points of order $p$.) An abelian variety $A$ of dimension $g$ over $\overline{\F}_q$ is \emph{ordinary} if $\#A(\overline{\F}_q)[p] = p^g$. 
Note that for dimension $g \geq 2$, there exist abelian varieties that are neither ordinary nor supersingular.

There is a natural isomorphism between the degree zero part of the Picard group $\Pic^0(C)$ of a hyperelliptic curve $C$ and its \emph{Jacobian} $\Jac_C$, which is an abelian variety into which the curve embeds (cf. \cite{FreyLangeCh4:2006}). For the remainder of this paper, we will identify the Picard group $\Pic^0(C)$ with $\Jac_C$.

\subsection{Arithmetic in the Jacobian}
\label{subsec: arithjac}

We will work in the Jacobian $\Jac_C$ of a hyperelliptic curve $\C$ of genus $g$, whose elements are equivalence classes of degree-zero divisors.  To do so, we choose a \emph{reduced} representative in each such divisor class.   A \emph{reduced} divisor is one of the form
\[
(P_1) + (P_2) + \cdots + (P_r) - r(P_\infty)
\]
where $r \leq g$, $P_\infty$ is the point at infinity on $\C$, $P_i \neq -P_j$ for distinct $i$ and $j$, and no $P_i$ satisfying $P_i = -P_i$ appears more than once.  Such a divisor is called \emph{semi-reduced} if the condition $r \leq g$ is omitted.  Each equivalence class contains exactly one reduced divisor.  For a divisor $D$ we will denote by $\rho(D)$ the reduced representative of its equivalence class.  The action of Galois commutes with $\rho$, i.e. $\rho(D^\sigma) = \rho(D)^\sigma$, since the property of being reduced is preserved by the action of Galois.

To encode the reduced divisor in a convenient way, we write $(u(x), v(x))$ where $u(x)$ is a monic polynomial whose roots are the $x$-coordinates $x_1, \ldots, x_r$ of the $r$ points
\[
P_1 = (x_1, y_1),\ \ldots, \ P_r = (x_r, y_r),
\]
and where $v(x_i) = y_i$ for $i=1,\ldots,r$.  This so-called \emph{Mumford representation} \cite{Mumford:83} is uniquely determined by and uniquely determines the divisor.  To find this representation, it suffices to find $u(x)$ and $v(x)$ satisfying the following conditions:
\begin{enumerate}
\item $u(x)$ is monic,
\item $\deg(v(x)) < \deg(u(x)) \leq g,$ and
\item $u(x) \mid F(x) - v(x)H(x) - v(x)^2,$
\end{enumerate}
where $F(x)$ and $H(x)$ are the polynomials defining the curve $\C$ (defined in Section \ref{hyperellipticcurvedefinition}).
When we add two reduced divisors $D_1$ and $D_2$ the result $D_1+D_2$ is not necessarily reduced.  Beginning with two reduced divisors in Mumford representation, the algorithm to obtain the Mumford representation of the reduction of their sum can be explained in terms of the polynomials involved in the Mumford representation, without recourse to the divisor representation.  This algorithm is originally due to Cantor \cite{Cantor:1987}, and in the form presented here to Koblitz \cite{Koblitz:89}.  The algorithm has two stages:  in the first, we find a semi-reduced divisor $D \sim D_1 + D_2$, and in the second stage, we reduce $D$.  Suppose that $D_i$ has Mumford representation $(u_i,v_i)$ for $i=1,2$.

\begin{stage}{1}
\begin{enumerate}
\item Find $d(x) = \gcd(u_1(x), u_2(x), v_1(x)+v_2(x)+H(x))$.  Finding this via the extended Euclidean algorithm gives $s_1(x)$, $s_2(x)$ and $s_3(x)$ such that
\[
d = s_1u_1 + s_2u_2 + s_3(v_1+v_2+H).
\]
\item Calculate the quantities
\[
u = u_1u_2/d^2, \ \ \mbox{and} \ \ v = s_1u_1v_2 + s_2u_2v_1 + s_3(v_1v_2 + F)/d \pmod{u(x)}.
\]
(It is easily verified that the fraction on the right is defined since $d(x)$ is a divisor of the numerator.)
\end{enumerate}
At this point, the result $(u, v)$ is a semi-reduced divisor linearly equivalent to $D_1+D_2$.  This stage corresponds to simply adding $D_1$ and $D_2$ and canceling any points with their negatives if applicable.  In fact, we obtain
\[
D' = D_1 + D_2 - \div(d).
\]
\end{stage}

\begin{stage}{2}
In this stage, if $\deg(u) > g$ we can replace $(u, v)$ with a divisor $(u', v')$ satisfying $\deg(u') < \deg(u)$.  This replacement is as follows. Set
\[
u' = (F - vH - v^2) / u, \ \ \mbox{and} \ \ v' = -H - v \pmod{u'}.
\]
This stage corresponds to simplifying the divisor using the geometric group law nicely described for genus $2$ by Lauter \cite{Lauter:2003}.  At each application of this loop to a divisor $D_3$, we obtain a divisor $D''$ satisfying\footnote{In general, $u'$ is a product of lines $L_i$ whose divisors are $(P_i) + (-P_i) - 2(P_\infty)$ for $i=1,\ldots, r$ and $\div(F-vH - v^2)$ is the sum of the intersection points of $\C$ and a unique curve intersecting $\C$ at $3g$ points including $P_1,\ldots,P_r$.}
\[
D'' = D_3 - \div( (F-vH - v^2) / u' ).
\]
Applying this loop finitely many times, beginning with the result $D'$ of stage one, we eventually obtain a reduced divisor $D$ linearly equivalent to $D_1 + D_2$. 
\end{stage}

This algorithm has been optimized to avoid the use of the extended Euclidean algorithm and in this form it is much more efficient \cite{GHV:2007}.  An enhanced version of Cantor's Algorithm is given as Algorithm \ref{algorithm: millerstep} in this paper; see Section \ref{subsec: miller}.  If steps 5 and 8 through 13 are removed from Algorithm \ref{algorithm: millerstep} one has the Cantor's Algorithm discussed here.

%*******************************************************

\section{Hyperelliptic Curves of Low Embedding Degree}
\label{s:friendly}

In this section we discuss hyperelliptic curves suitable for pairing-based cryptographic systems. The Jacobian varieties of such curves must have computable pairings, and computationally infeasible discrete logarithm problems. Specifically, we require low embedding degrees and large prime-order subgroups. 

\subsection{Embedding degree and $\rho$-value.} 
Let $r$ be a prime.
Let $C$ be a hyperelliptic curve over $\F_q$ of genus $g$ with Jacobian variety $\jac_C(\F_q)$ such that $r\mid \#\jac_C(\F_q)$ and $\gcd(r,q) =1$. The {\em embedding degree} of $\jac_C$ with respect to $r$ is the smallest integer $k$ such that $r\mid (q^k-1)$. Equivalently, the embedding degree of $\jac_C$ with respect to $r$ is the smallest integer $k$ such that $\F^\ast_{q^k}$ contains the group of $r^{\th}$ roots of unity $\mu_r$. If $\Jac_C$ has embedding degree $k$ with respect to $r$, then a pairing on $C$, such as the Weil pairing $e_r: \Jac_C(\F_q)[r]\times \Jac_C(\F_q)[r] \to \mu_r$, ``embeds'' $\Jac_C(\F_q)[r]$ (and any discrete logarithm problem in $\Jac_C(\F_q)[r]$) into $\F^\ast_{q^k}$, and $\F_{q^k}$ is the smallest-degree extension of $\F_q$ with this property; whence the name ``embedding degree''.
Hitt \cite{Hitt:2007} shows that if $q=p^m$ with $m>1$, then $\Jac_C(\F_q)[r]$ may be embedded into a smaller field which is not an extension of $\F_q$ but only an  extension of $\F_p$. The smallest such field is the so-called {\em minimal embedding field}, which is $\mathbb{F}_{p^{\ord_r p}}$.

We occasionally speak of the embedding degree of the hyperelliptic curve $C$, in which case we mean the embedding degree of its Jacobian. 

Another important parameter is the {\em $\rho$-value}, which for a Jacobian variety of dimension $g$ we define as $\rho = g\log q/\log r$. Since $\#\Jac_C(\F_q)= q^g  + O(q^{g-1/2})$, the $\rho$-value  measures the ratio of the bit-sizes of $\#\jac_C(\F_q)$ and the subgroup order $r$. Jacobian varieties with a prime number of points have the smallest $\rho$-values: $\rho\approx 1$.  
We call a hyperelliptic curve, and its Jacobian variety, {\em pairing-friendly} if the Jacobian variety  has small embedding degree and a large prime-order subgroup. In practice, we want  $k\le 60$ and  $r> 2^{160}$.

Since the embedding degree $k$ is the order of $q$ in the multiplicative group $(\Z/r\Z)^\ast$, and typically elements in $(\Z/r\Z)^\ast$ have large order, we expect that for a random Jacobian over $\F_q$ with order-$r$ subgroup,  the embedding degree is approximately of the same size as $r$. (This reasoning has been made more precise for elliptic curves, by Balasubramanian and Koblitz \cite{BalasubramanianKoblitz:1998} and Luca, Mireles and Shparlinski \cite{LMS:2004}.) With $r>2^{160}$, this means that evaluating a pairing for a {\em random} hyperelliptic curve becomes a computationally infeasible task. Just as in the case of elliptic curves, pairing-friendly hyperelliptic curves are rare and require special constructions. 

\subsection{Embedding degrees required for various security levels.}
For cryptographic applications,  the discrete logarithm problems in $\jac_C(\F_q)$ and in the multiplicative group $\F_{q^k}^\ast$ must both be computationally infeasible.
For Jacobian varieties of hyperelliptic curves of genus $2$ 
 the best known discrete logarithm (DL) algorithm is the parallelized Pollard rho algorithm 
\cite{OorschotWiener:99, Pollard:78}, which has running time $O(\sqrt{r})$ where $r$ is the size of the largest prime-order subgroup of $\jac_C(\F_q)$.
For Jacobian varieties of dimensions $3$ and $4$ there exist index calculus algorithms of complexities $O(q^{4/3+\varepsilon})= O(|\jac_C|^{4/9+\varepsilon})$ and $O(q^{3/2+\varepsilon}) = O(|\jac_C|^{3/8+\varepsilon})$, respectively \cite{GTTD:2007}.
How this compares to the parallelized Pollard rho algorithm 
depends on the relative size of the subgroup order $r$ -- more precisely, only if $\rho < 9/8$ (genus $3$ case) or $\rho < 4/3$ (genus $4$ case) will the index calculus approach be superior to Pollard rho. 

In any case, the best DL algorithms for genus $2,3,$ and $4$ are of exponential running time. On the other hand, the best algorithm for DL computation in finite fields is the index calculus attack (e.g., \cite{Odlyzko:1985}) which has running time subexponential in the field size.  Thus to achieve the same level of security in both groups, the size $q^k$ of the extension field must be significantly larger than $r$. Table \ref{tab:keys} shows sample subgroup sizes, extension field sizes, and embedding degrees with which to achieve common levels of security, for various cases $r\approx q^{g/\rho}$. 
The listed sizes for the prime-order subgroups and the extension fields (of large characteristic) follow the recommendations by NIST \cite[Table 2]{NIST:2006}. 

\begin{table}[ht]
\caption{Embedding degrees for hyperelliptic curves of genus $g=2$ required to obtain commonly desired levels of security.}%in bits.}
\label{tab:keys}
\begin{tabular}{|c|c|c|c|c|c|c|c|c|c|} \hline
Security & Subgroup  & Extension field & 
\multicolumn{6}{c|}{Embedding degree $(k)$}  \\ 
level (bits) & size $(r)$ & size $(q^k)$ & $\rho \approx 1$ & $\rho \approx 2$ & $\rho\approx 3$ & $\rho\approx 4$ & $\rho\approx 6$ & $\rho\approx 8$\\ \hline
 80 & 160  &  1024	&  $6g$	 	 & $3g$  & $2g$ & $1.5g$ & $g$  &  $0.8g$ \\  
 112 & 224 & 2048	&  $10g$	 & $5g$  & $3.3g$ & $2.5g$ & $1.6 g$ &  $1.3g$ \\  
 128 & 256 & 3072	&  $12g$	 & $6g$  & $4g$ & $3g$ & $2g$  & $1.5g$ \\ 
 192 & 384 & 7680	&  $20g$	 & $10g$ & $6.6g$ & $5g$ & $3.3g$ & $2.5g$ \\
 256 & 512 & 15360      & $30g$          & $15g$ & $10g$ & $7.5g$ & $5g$ &  $3.8g$ \\ \hline 
\end{tabular}
\vspace{10pt}
\end{table}
% binary field sizes for these 5 security levels (following Lenstra's 2001 asiacrypt paper): 1500, 3000, 4500, 11250, 22500

While Table \ref{tab:keys} as such is for genus $2$ only, it can easily be adapted to the cases of genus $3$ and $4$: Only in the case that the Jacobian has almost prime order ($\rho \approx 1$) we need to compensate for the aforementioned  index-calculus algorithms in $\jac_C$. For this, if $g=3$, multiply the second column entries by $9/8$ and the fourth column entries by $8/9$; if $g=4$ multiply the second column entries by $4/3$ and the fourth column entries by $3/4$.

\subsection{Ordinary hyperelliptic curves of low embedding degree} 
\label{ss:ordinary}
While there are numerous constructions for pairing-friendly elliptic curves -- see e.g.\ the survey by Freeman, Scott and Teske \cite{FST:2009} --  there are not nearly as many constructions for hyperelliptic curves of low embedding degree and large prime-order subgroup. In this section, we discuss the case of ordinary Jacobians; see Section \ref{ss:supersingular} for the supersingular case. We keep the discussion result-oriented, and refer the reader to the corresponding original papers for details on the specific constructions and the theory underneath. 

Galbraith, McKee and Valen\c ca \cite{GalbraithMcKeeValenca:2005} argue that heuristically, for any fixed embedding degree $k$ with $\varphi(k) \ge 4$ ($\varphi(k) = $ the Euler phi-function) and for any bound $M$ on the field size $q$, there exist about as many genus $2$ curves over $\F_q$ of embedding degree $k$ (any $\rho$-value) as there exist elliptic curves over $\F_q$ of embedding degree $k$, namely $\Theta(M^{1/2}/\log M)$. For embedding degrees $k=5,10$, they identify several quadratic polynomials $q(x)$ parameterizing field sizes such that genus $2$ curves over $\F_{q(x)}$ exist with embedding degree $k$ (any $\rho$-value). (They also show that for $k=8,12$, such quadratic polynomials $q(x)$ do not exist.) 

Freeman \cite{Freeman:2007} was the first to actually construct ordinary genus 2 curves of low embedding degree. His construction is based on the Cocks-Pinch method  \cite{CocksPinch:2001}\cite[Theorem 4.1]{FST:2009}, which produces pairing-friendly {\em elliptic} curves over prime fields of any prescribed embedding degree and with $\rho \approx 2$. In the genus-2 case, Freeman obtains curves over prime fields $\F_q$ of any prescribed embedding degree $k$ and $\rho$-value $8$, that is, $r\approx q^{1/4}$ (where $r$ denotes the prime subgroup order of the Jacobian).

Freeman  \cite[Proposition 2.3]{Freeman:2007} further shows that the resulting Jacobian varieties have the property that $\Jac_C(\F_{q^k})$ always contains two linearly independent $r$-torsion points. For an elliptic curve $E/\F_q$, the corresponding result implies that the entire $r$-torsion group is contained in $E(\F_{q^k})$, but this is not necessarily the case for higher dimensional abelian varieties. This phenomenon gives rise to the notion of the {\em full} embedding degree, which is the smallest integer $k$ such that all $r$-torsion points of $\Jac_{\C}$ are defined over $\F_{q^k}$. Freeman \cite[Algorithm 5.1]{Freeman:2007} gives a construction of genus 2 curves of prescribed full embedding degree $k$ (necessarily even), which may be useful in cryptographic applications that require more than two linearly independent $r$-torsion points (see Section \ref{ss:highdimTorsion}). Again, this construction yields curves with $\rho$-value $8$. 

Note that an essential part of either construction \cite{Freeman:2007} is the use of the complex multiplication (CM) method to compute the actual curve. In genus 2, this includes computation of the Igusa class polynomials (e.g., \cite{Weng:2003}) of 
the CM  field $K = \rm{End}(\Jac_{\C}) \otimes \Q$, which is currently feasible for CM fields $K$ with class numbers less than $100$ \cite{KohelDataBase}. (Here, $\rm{End}(\Jac_{\C})$ denotes the set of all endomorphisms of $\Jac_C$ defined over $\F_q$.)

Freeman, Stevenhagen and Streng  \cite[Algorithm 2.12]{FSS:2008} present a generalization of the Cocks-Pinch method, which, when coupled with complex multiplication methods, produces pairing-friendly abelian varieties over prime fields,  of dimension $g$ with $\rho$-values $\approx 2g^2$. This algorithm works for any prescribed embedding degree $k$, and applies to arbitrary genus $g\ge 2$. (However note that complex multiplication methods are available for special CM fields only if $g=3$, and are completely undeveloped for $g\ge 4$.) In addition to explicit genus $2$ examples with $\rho\approx 8$, a cryptographically interesting example is given for genus $3$ ($k=17$ and $\rho \approx 17.95$).

In the case of pairing-friendly {\em elliptic} curves, the method by Brezing and Weng \cite{BrezingWeng:2005} is a generalization of the Cocks-Pinch method \cite{CocksPinch:2001} and produces elliptic curves over prime fields with $1<\rho<2$ for many embedding degrees. Freeman \cite[Algorithm 3.8]{Freeman:2008} combines the Brezing-Weng approach
with the method from Freeman, Stevenhagen and Streng \cite{FSS:2008} to construct so-called {\em families} of abelian varieties over prime fields with $\rho$-values strictly less than $2g^2$.
An explicit construction for genus $2$, embedding degree $k=5$ and $\rho = 4$ is given -- note that an instantiation with a $224$-bit prime subgroup order $r$ would exactly meet the $112$-bit security level requirements (cf.\ Table \ref{tab:keys}).  Other examples (for genus $2$) include: $k=6$, $\rho = 7.5$; $k=8$, $\rho = 7.5$, and $k=10$, $\rho=6$ (able to exactly meet the $256$-bit security level requirements) \cite{Freeman:2008, FreemanDataBase}.
In the case of genus $3$, a construction yielding $k=7$ and $\rho=12$ is obtained. 

All constructions mentioned so far in this section produce absolutely simple Jacobians. When  considering simple abelian varieties $A$ that are
isogenous over some extension field $\F_{q^d}$ ($q$ a prime) to a product of two elliptic curves, smaller $\rho$-values have been obtained:\\ 
Kawazoe and Takahashi \cite{KawazoeTakahashi:2008}  
specialize to hyperelliptic curves with curve equation $y^2 = x^5+ax$  
over a prime field $\F_q$. For the cardinalities of the Jacobians of such curves, closed formulae exist. These
formulae are exploited in adaptations of the Cocks-Pinch method (producing Jacobians with $\rho$-values around $4$), and  
Brezing-Weng-type methods (for embedding degree divisible by $8$, producing Jacobian varieties with  
$3<\rho<4$). The Jacobians split over $\F_{q^d}$, $d\in \{2,4\}$.\\
Satoh  \cite{Satoh:2009} considers 
hyperelliptic curves $C$ of the form $y^2=x^5+ax^3+bx$ over $\F_q$, such that $\Jac_C$ splits over $\F_{q^2}$. This construction works for many embedding degrees and produces $\rho$-values $< 4$.\\
More generally, Freeman and Satoh \cite{FreemanSatoh} show that if $E$ is
defined over $\F_q$, and $A$ is an abelian variety isogenous over $\F_{q^d}$ to $E\times E$, then $A$ is isogenous over $\F_q$ to a primitive subvariety of the Weil restriction of $E$ from $\F_{q^d}$ to
$\F_q$. Thus, pairing-friendly abelian varieties of this type can be built from elliptic curves $E/\F_q$ that are not
pairing-friendly over $\F_q$, but are pairing-friendly when base-extended
to $\F_{q^d}$. The elliptic curves can be constructed via Cocks-Pinch or Brezing-Weng type methods. The
generic $\rho$-value for Jacobians of genus $2$ produced in this manner is $4$.
With the Brezing-Weng method, $\rho$-values between $2$ and $4$ can be obtained. This approach not only contains the constructions by Kawazoe and Takahashi \cite{KawazoeTakahashi:2008} and Satoh \cite{Satoh:2009} but also produces the lowest ever recorded $\rho$-values for ordinary genus $2$ curves. 
Explicit examples of cryptographically interesting genus $2$ curves are given, such as a $k=9$, $\rho\approx 8/3$ curve and a $k=27$, $\rho \approx 20/9$ curve.

In conclusion, to date, the best we can achieve for pairing-friendly ordinary genus $2$ curves with arbitrary prescribed embedding degree $k$ is a $\rho$-value of $4$; and $\rho\approx 8$ if one insists on absolutely simple Jacobians. (Although to date, there is no apparent reason why Jacobians that split over small-degree extensions should be more vulnerable to DL attacks than the absolutely simple ones.)
We have no constructions of ordinary hyperelliptic curves of genus $g\ge 2$  with $\rho$-values less than $2$. In particular, we have no constructions of higher-dimensional pairing-friendly ordinary Jacobian varieties with a prime number of points. 
This is in sharp contrast to the elliptic case, where $\rho\approx 2$ can be achieved for any prescribed embedding degree, $1 < \rho <2$ for selected embedding degrees,  and constructions for prime-order elliptic curves exist for embedding degrees $k=3,4,6,10,$ and $12$ (cf.\ \cite{FST:2009}).

\subsection{Supersingular curves}
\label{ss:supersingular}
Supersingular hyperelliptic curves over $\F_q$ are always
pairing-friendly. In fact, Galbraith \cite{Galbraith:2001} shows that
there exists a constant $k(g)$ such that the embedding degree of any
supersingular abelian variety of dimension $g$ over any finite field
$\F_q$ is bounded by $k(g)$. Rubin and Silverberg
\cite{RubinSilverberg:2002} prove that for simple supersingular
abelian varieties, for $g\le 6$ we have $k(g) \le 7.5g$.

Specifically, for dimension $g=2$, the embedding degree is bounded by $12$, where $k=12$ can only happen if $\F_q$ is a binary field $\F_{2^m}$ with $m$ odd. If $q$ is a square, or if $q=p^m$ with $m$ odd and $p\neq 2,3$, then the largest embedding degree is $k=6$. If $\F_q = \F_{3^m}$ with $m$ odd, the embedding degree is always bounded by $4$. 
(In the case of dimension $g=3$, the embedding degree is bounded by $18$, and the bound for the dimension $4$ case is $30$. In both cases, this bound is achieved only in characteristic three. Over prime fields $\F_p$ with $p\ge 11$, there are no simple supersingular abelian varieties of dimension $g=3$, while the largest embedding degree for dimension $g=4$ is $k=12$.)

As Rubin and Silverberg show \cite[Corollaries 13,14]{RubinSilverberg:2002}, not all embedding degrees below these bounds are possible. For example, in the dimension  $2$ case and if $q=p^m$ with $m$ odd, then for $p=2$ we have $k\in \{1,3,6,12\}$;  if $p=3$, we have $k\in \{1,3,4\}$; if $p=5$ we have $k\in \{1,3,4,5,6\}$ and if $p\ge 7$ we have $k\in\{1,3,4,6\}$. 

Cryptographically interesting supersingular hyperelliptic curves can
be explicitly constructed. For example, Galbraith et al.\
\cite{GPRS:2009} give curve equations for various field
characteristics that yield simple supersingular Jacobians of dimension
$g=2$ and of embedding degrees $k\in \{4,5,6,12\}$. By carefully
choosing the underlying fields, $\rho$-values close to $1$ can be
readily obtained.

\subsection{Supersingular versus ordinary hyperelliptic curves}
While the embedding degrees of supersingular abelian varieties are limited to a few, small values, their advantage is that they can achieve $\rho$-values significantly smaller than their ordinary counterparts. For example, let us consider the $112$-bit security level (cf.\ Table \ref{tab:keys}). One could use the construction by Freeman and Satoh \cite{FreemanSatoh} of an ordinary absolutely simple hyperelliptic Jacobian of dimension $2$, with embedding degree $k=6$ and $\rho$-value $2.976$, with a $230$-bit prime-order subgroup, working over a finite field $\F_q$ with $342$-bit $q$. Alternatively, one could use the embedding-degree $12$ supersingular curve $y^2+y =x^5+x^3+b$ ($b\in\{0,1\}$) over $\F_{2^m}$ with $m\ge 250$ chosen such that its Jacobian contains a subgroup of prime order $r>2^{224}$. 
(Note that Coppersmith' algorithm \cite{Coppersmith:1984} for DL computation in finite fields of small characteristic requires to embed the Jacobian into a $3000$-bit binary field $\F_{2^{12m}}$, to obtain roughly the same level of security provided by a $2048$-bit field $\F_{q^{12}}$ with $q$ large, cf.\ \cite{Lenstra:2001}.)
If $m$ is chosen smaller than $342$, this would result in bandwidth advantages for the supersingular Jacobian, given that in cryptographic applications the values that are transmitted are elements in $\Jac_C(\F_q)$. 
However, already at the $128$-bit security level the advantage of supersingular curves disappears, in the light of the recent work by Freeman and Satoh \cite{FreemanSatoh}:
this security level can be achieved with $256$-bit prime-order subgroups  either of an ordinary Jacobian over a $341$-bit $\F_q$, with $k=9$ and $\rho = 8/3$, or of a supersingular Jacobian over $\F_{2^m}$ with $m\ge 375$, of embedding degree $12$ (again, $m$ is chosen in response to Coppersmith' DL algorithm \cite{Coppersmith:1984}: a $4500$-bit binary field roughly provides the same security as a $3072$-bit field of large characteristic). At high security
levels ordinary curves are definitely preferable. For example, at the $256$-bit level, a genus $2$ curve with embedding degree $k=27$ and  (optimal to date) $\rho$-value of $20/9$ (cf.\ \cite{FreemanSatoh})  requires a $568$-bit field, while a binary  supersingular curve of embedding degree $12$ requires a $1875$-bit field.

%*******************************************************

\section{Pairings for Hyperelliptic Curves}
\label{s:pairings}

In this section, we give an overview of the different pairings on hyperelliptic curves, as well as introduce the
more general framework of {\it HV pairings} which unify the recent variations on the Ate pairing. In particular, we present a direct
proof of bilinearity and non-degeneracy for these pairings and describe how the Ate$_i$ and R-ate pairings fit into the framework.

We begin by introducing  the historically most important pairings for hyperelliptic curves, the Tate-Lichtenbaum and Weil pairings. In what follows, let $r$ be a positive integer and assume that $\C$ is defined over a finite field $\F_q$.  Suppose that $K = \F_{q^k}$ is an extension of $\F_q$ such that $r \mid (q^k-1)$. Throughout the section, we will use $D$ to mean both a divisor and the  divisor class represented by $D$.

For a positive integer $s$,  a {\it Miller function} $f_{s,D}$ is a function with divisor
\[
(f_{s,D}) = sD - \rho(sD),
\]
uniquely defined up to scalar multiplication by elements of $K^*$. 
The {\it Miller loop length} of such a function is $\log_2 s$ and measures how quickly the function can be evaluated via
Miller's algorithm (see Algorithm \ref{algorithm: miller}).  The benefit of recent variations on the Tate-Lichtenbaum pairing is a reduction in  Miller loop length, which is sometimes accomplished by combining several Miller functions (see Section \ref{s:fast}).

\subsection{Tate-Lichtenbaum pairing}

For $D_1 \in \Jac_C(K)[r]$, the divisor $rD_1$ is linearly equivalent to zero, hence there is some function whose divisor is $rD_1$, namely the Miller function $f_{r,D_1}$ defined above.  Let $D_2$ be a divisor class, with representative $D_2 = \sum_P n_P(P)$ disjoint from $D_1$.
 We define a pairing called the \emph{Tate-Lichtenbaum pairing} as follows
\begin{eqnarray*}
\tau : \Jac_C(K)[r] \times \Jac_C(K) / r\Jac_C(K) &\rightarrow& K^*/(K^*)^r \\
(D_1,D_2) &\mapsto& f_{r,D_1}(D_2) = \prod_P f_{r,D_1}(P)^{n_P}.
\end{eqnarray*}
This pairing is bilinear, non-degenerate and the result is independent of the choice of representatives of the divisor classes.

\subsection{The Weil pairing}

For $D_1, D_2 \in \Jac_C(\bar K)[r]$, the \emph{Weil pairing} is given by
\begin{eqnarray*}
e_r:  \Jac_C(\bar K)[r] \times \Jac_C(\bar K)[r] &\rightarrow& \mu_r \\
(D_1,D_2) &\mapsto& \tau(D_1,D_2) \tau(D_2,D_1)^{-1}
\end{eqnarray*}
which can be computed via two Tate-Lichtenbaum pairings. It is bilinear, alternating, and non-degenerate.

\subsection{The modified Tate-Lichtenbaum pairing}

If $\Jac_C(K)$ contains no elements of order $r^2$, then there is an isomorphism
\[
\Jac_C(K)[r] \cong \Jac_C(K) / r \Jac_C(K).
\]
Under this identification, we define the \emph{modified (or reduced) Tate-Lichtenbaum pairing} to be
\begin{eqnarray*}
t : \Jac_C(K)[r] \times \Jac_C(K)[r] &\rightarrow& \mu_r \\
(D_1,D_2) &\mapsto& \tau(D_1,D_2)^{(q^k-1)/r}.
\end{eqnarray*}

Since elements of $K^*$ have order dividing $q^k-1$ and $r \mid (q^k-1)$, the $r^{\th}$ powers which are the quotients of distinct representatives of the coset of $\tau(D_1,D_2)$ are removed by this {\it final exponentiation}, leaving a unique result lying in $\mu_r \subset K$.

Other powers of the Tate-Lichtenbaum pairing can also give non-degenerate bilinear pairings into $\mu_r$ which may yield shorter Miller loops (for example, with the use of efficiently computable automorphisms of $\C$ \cite{FGJ:2008}; see Section \ref{ss:automorphisms_faster}).

\subsection{Hyperelliptic Ate pairing}
\label{s:hyperate}

More generally, a \emph{bilinear pairing} is a map
\[
e : \G_1 \times \G_2 \rightarrow \G_3
\]
where $\G_i$ are abelian groups, in additive notation,  and $\G_3$ is a cyclic group, written multiplicatively, and for all $p_1,p_2 \in \G_1$, $q_1,q_2 \in \G_2$, we have
\begin{align*}
e(p_1 + p_2, q_1) &= e(p_1, q_1)e(p_2, q_1), \\
e(p_1, q_1 + q_2) &= e(p_1, q_1)e(p_1, q_2).
\end{align*}

Let $r$ be a prime dividing $\#\Jac_C(\F_q)$ and let $k$ be the embedding degree of $\Jac_C(\F_q)$ with respect to $r$. We are interested in pairings where $\G_1$ and $\G_2$ are subgroups of $\Jac_C(K)$, 
where $K = \F_{q^k}$.  In particular, a number of more convenient and faster pairings are known when
\begin{equation}
\label{def:G1G2}
\begin{array}{ll}
\displaystyle \G_1 &= \Jac_C(K)[r] \cap \ker( \pi - [1] ), \\
\displaystyle \G_2 &= \Jac_C(K)[r] \cap \ker( \pi - [q] ),
\end{array}
\end{equation}
where $\pi$ is the $q^{\th}$ power Frobenius automorphism. Since $r$ divides $\#\Jac_C(\F_q)$, the group $\G_1$, being the eigenspace of $1$, is at least $1$-dimensional over $\Z/r\Z$. Since the eigenvalues of the Frobenius come in pairs $(\lambda, q/\lambda)$ \cite[\S5.2.3]{FreyLangeCh5:2006}, $q$ is also an eigenvalue of $\pi$ on $\Jac_C[r]$, and thus there exists a divisor $D$ such that  $\pi(D) = qD$. This implies that  $\pi^k D = q^k D = D$, since $r|(q^k -1)$ and $rD = 0$. Consequently, $D \in \Jac_C(\F_{q^k})$, and the group $\G_2$ is also at least $1$-dimensional over $\Z/r\Z$.
If $k>1$, then $\G_1 \neq \G_2$ and  $\G_1 \times \G_2 \subset \Jac_C(\F_{q^k})[r]$ is at least 2-dimensional over $\Z/r\Z$. (Recall that for genus $g$, the group $\Jac_C(\overline K)[r]$ is $2g$-dimensional over $\Z/r\Z$.)

In the remainder of this section, $\G_1$ and $\G_2$ always denote the groups defined in \eqref{def:G1G2}. 

The most basic pairing defined for divisors in $\G_1$, $\G_2$ is the {\it hyperelliptic Ate pairing} \cite{GHOTV:2007}:
\begin{eqnarray*}
a:\G_2 \times \G_1 &\rightarrow& \mu_r \\
(D_2,D_1) &\mapsto& f_{q, \rho(D_2)}(D_1),
\end{eqnarray*}
where $\rho(D_2)$ is the reduced divisor class representative. 
Since the Frobenius $\pi$ acts as $[q]$ on $D_2$, we have
$f_{q,\rho(D_2)}(D_1) \in \mu_r$ and no final exponentiation is required \cite[Lemma 2]{GHOTV:2007}.  This is different from the elliptic Ate pairing \cite{HSV:2006}, where a final exponentiation is always required.
Another important difference of the hyperelliptic Ate pairing is that to obtain a well-defined value, one \textit{must} use the reduced divisor $\rho(D_2)$, not simply any representative of the class $D_2$. 
The Miller loop length for the hyperelliptic Ate pairing is $\log_2 q$, in contrast to the elliptic case where
the Miller loop length is $\log_2 (t-1)$ with $t$ the trace of Frobenius.

\subsection{The Hess-Vercauteren (HV) framework for pairings on Frobenius eigenspaces}
\label{ss:framework}

Since 2007, several variations of the Ate pairing have been proposed for elliptic and hyperelliptic curves, exploiting the fact that products and ratios of bilinear, non-degenerate pairings on $\G_2 \times \G_1$ are also bilinear pairings, but not necessarily non-degenerate \cite{ZZH:2008-2}. The key is to find combinations of pairings which are both non-degenerate and computable using shorter Miller loops. Following  the work of Hess \cite{Hess:2008} and Vercauteren \cite{Vercauteren:2008} in the elliptic curve case, we unify these various pairings on $\G_2 \times \G_1$ in a more general framework, which we call {\it HV pairings}. The main benefit of this framework is that the criteria for non-degeneracy are more straightforward to verify, giving a direct way to create new pairings. Further investigation of this framework and 
possible extensions seems likely to be fruitful (see Section \ref{ss:optimalpairings} and $(1)$ in Section \ref{ss:other}).

Let $D$ be any divisor in $\Jac_C(K)[r]$, and $s$ an integer.  Recall that any divisor $D$ is equivalent to a unique reduced divisor which we denote $\rho(D)$.
Let $h(x) \in \mathbb{Z}[x]$ be a polynomial of the form $h(x) = \sum_{i=0}^{n} h_i x^i$ satisfying $h(s) \equiv 0 \pmod{r}$.  Define a \emph{generalized Miller function} $f_{s,h,D}$ to be any function with divisor
\begin{equation}
\label{generalisedMillerfunction}
\sum_{i=0}^{n} h_i \rho(s^iD).
\end{equation}
To see that this divisor is principal, consider  the principal divisor
\[
\sum_{i=0}^{n} h_i (s^iD - \rho(s^iD)),
\]
which differs by $(\sum_{i=0}^{n} h_i s^i)D$ from (\ref{generalisedMillerfunction}). Since $h(s) \equiv 0 \pmod{r}$, this  is an integer multiple of $rD$, which is linearly equivalent to zero by assumption, and thus the divisor (\ref{generalisedMillerfunction}) is principal.  As with the standard Miller function, the function $f_{s,h,D}$ is only defined up to scalar multiples. Also, we note that the Miller function  $f_{r,D}$ for the Tate-Lichtenbaum pairing is equal to $f_{s,h,D}$ for the constant function $h(x) = r$ and arbitrary integer $s$.

\begin{theorem}
\label{HVtheorem}
Let $s \equiv q^j \pmod{r}$ for some $j \in \Z$.  Let $h(x) \in \mathbb{Z}[x]$ with $h(s) \equiv 0 \pmod{r}$.  Then
\begin{eqnarray*}
a_{s,h}:\mathbb{G}_2\times\mathbb{G}_1 &\rightarrow& \mu_r \\
(D_2,D_1) &\mapsto& f_{s,h,D_2}(D_1)^{(q^k-1)/r}
\end{eqnarray*}
is a  bilinear pairing satisfying
\[
a_{s,h}(D_2,D_1) = t(D_2,D_1)^{h(s)/r} \text{ and }  a_{s,h}(D_2, D_1) = a(D_2, D_1)^{kq^{k-1}h(s)/r}
\]
where $t$ is the modified Tate-Lichtenbaum pairing and $a$ is the hyperelliptic Ate pairing.
The pairing $a_{s,h}$ is non-degenerate if and only if $h(s) \not \equiv 0 \pmod{r^2}$.
\end{theorem}

\begin{remark}
We note that since  $k$ is the embedding degree of $\jac_C(\F_q)$ with respect to $r$, in Theorem \ref{HVtheorem} $s$ will be a $k^{\th}$ root of unity modulo $r$ since $q$ is a primitive $k^{\th}$ root. In Hess's framework, there is the additional  condition that $s$ be a primitive $k^{\th}$ root of unity modulo $r^2$. This requirement is necessary to show the existence of pairings such that the function $f_{s,h,D}$ is of ``lowest degree'' (see \cite[\S 3]{Hess:2008}), 
but is not required for the result above.
\end{remark}

\begin{proof}
First we show that the pairing is well-defined on divisor classes.  Suppose that $D_2' \sim D_2$.  Then
\[
\div(f_{s,h,D_2'}/f_{s,h,D_2}) = \sum_{i=1}^n h_i (\rho(s^iD_2') - \rho(s^iD_2)) = \emptyset.
\]
This demonstrates well-definition in the factor $\mathbb{G}_2$.  For the factor $\mathbb{G}_1$, it suffices to show that the pairing is trivial under the hypothesis that $D_1$ is a principal divisor.  Suppose $D_1 = \div(g)$.  For any $D_2 \in \mathbb{G}_2$, by the hypothesis that $s \equiv q^j \mod r$ and $\rho(rD_2)= \emptyset$, it is the case that
\[
\rho(s^iD_2) = \rho(q^{ij}D_2) = \rho(\pi^{ij}D_2) = \pi^{ij}\rho(D_2)= q^{ij}\rho(D_2) = s^i\rho(D_2) + rD'
\]
for some divisor $D'$ defined over $\mathbb{F}_q$.  Therefore,
\[
\sum_{i=1}^n h_i \rho(s^iD_2) = \sum_{i=1}^n h_i s^i\rho(D_2) + rD''
\]
for some $D''$ defined over $\mathbb{F}_q$.  Then by the hypothesis, this expression is an $r$-th multiple of another divisor $D'''$ defined over $\mathbb{F}_q$.  By Weil reciprocity,
\[
f_{s,h,D_2}(D_1)^{(q^k-1)/r} = g(rD''')^{(q^k-1)/r} = g(D''')^{q^k-1} = 1,
\]
as required.

We show bilinearity and non-degeneracy directly, in contrast to Hess's more general approach in the elliptic curve case \cite[Theorem 1]{Hess:2008}.

Let $s = q^j + \ell r$, for $j, \ell \in \Z$. Linearity in the second coordinate follows from the definition of evaluation of a function on a divisor. To show linearity in the first coordinate, 
let $D_2, D_3 \in \G_2$ and $D_1 \in \G_1$ be non-trivial reduced divisors. Then
\[
(f_{s,h, D_2 + D_3}) = \sum_{i = 0}^{n} h_i \rho(s^iD_2 + s^iD_3) = \sum_{i=0}^{n} h_i \rho(s^iD_2) + \sum_{i=0}^{n} h_i \rho(s^iD_3) + \sum_{i = 0}^{n} h_i (g_i)
\]
where
\[
(g_i) = \rho(s^iD_2 + s^iD_3) - \rho(s^iD_2) - \rho(s^iD_3).
\]
Since $rD_2 \sim 0$, $rD_3 \sim 0$ and $s = q^j + \ell r$, the function $g_i$ has divisor
\[
(g_i) = \rho(q^{ij}D_2 + q^{ij}D_3) - \rho(q^{ij}D_2) - \rho(q^{ij}D_3).
\]
Since $D_2, D_3 \in \G_2$, the $q$-eigenspace of the Frobenius $\pi$, and since $\rho$ commutes with $\pi$, we have
\[
(g_i) = \rho(D_2 + D_3)^{\pi^{ij}} - \rho(D_2)^{\pi^{ij}} - \rho(D_3)^{\pi^{ij}}.
\]
Then $(g_i) = (m)^{\pi^{ij}}$ where $m$ is the function  with divisor
\[
(m) = \rho(D_2 + D_3) - \rho(D_2) - \rho(D_3).
\]
As $f_{s,h,D_2 + D_3}$ is evaluated at the divisor $D_1 \in \G_1$, which is fixed by $\pi$, the value $g_i(D_1)$ equals $m(D_1)^{\pi^{ij}} = m(D_1)^{q^{ij}} $.
Thus,
\[
\prod_{i=0}^{n} g_i(D_1)^{h_i} = \prod_{i=0}^{n} m(D_1)^{h_iq^{ij}} = m(D_1)^{\sum_{i=0}^{n} h_iq^{ij}} = m(D_1)^{h(q^{j})}.
\]
Using the fact that $s = q^j + \ell r$ and $h(s) \equiv 0 \pmod{r}$, we see that this value is eliminated by the final exponentiation of $(q^k-1)/r$.
Since
\[
f_{s,h,D_2 + D_3}(D_1)  =  f_{s,h,D_2}(D_1)f_{s,h,D_3}(D_1) \prod_{i=0}^{n} g_i(D_1)^{h_i} ,
\]
the pairing $a_{s,h}$ is linear with respect to the first coordinate.

We now show that
\[
a_{s,h}(D_2,D_1) = t(D_2,D_1)^{h(s)/r}
\]
using a similar argument.
On the right, we have
\[
t(D_2,D_1)^{h(s)/r} = \left( f_{r,D_2}(D_1)^{(q^k-1)/r} \right) ^{h(s)/r}.
\]
Since $D_2 \in \G_2$, we have $\rho(rD_2) = 0$, thus
\[
(f_{r,D_2}^{\,h(s)/r}) = (h(s)/r) ( rD_2 - \rho(rD_2) ) = h(s)D_2 = \sum_{i = 0}^{n} h_i s^i D_2.
\]

On the left, we have
\[
a_{s,h}(D_2,D_1) = f_{s,h,D_2}(D_1)^{(q^k-1)/r}.
\]
where by definition
\[
(f_{s,h,D_2}) =  \sum_{i = 0}^{n} h_i  \rho(s^iD_2).
\]
We can rewrite this as
\[
(f_{s,h,D_2}) = \sum_{i = 0}^{n} h_i s^i D_2 -  \sum_{i = 0}^{n} h_i(g_i),
\]
where
\[
(g_i) = s^iD_2 -\rho(s^i D_2).
\]
Since we evaluate at  $D_1 \in \G_1$ fixed by $\pi$ and  $s = q^j + \ell r$ for some $\ell \in \Z$, the contribution of the function with divisor $(\sum_{i = 0}^{n} h_i(g_i))$ is eliminated by
 raising to the power $(q^k - 1)/r$. Furthermore, we may choose any functions $f_{r,D_2}$ and $ f_{s,h,D_2}$ with the above divisors, as any discrepancy from scalar multiples will be 
 canceled out when evaluating at the degree zero divisor $D_1$.
 Thus,  $a_{s,h}(D_2,D_1) = t(D_2,D_1)^{h(s)/r}$.

We have that $t$ is a non-degenerate pairing and $h(s) \equiv 0 \pmod{r}$. Therefore, by the relationship between $a_{s,h}$ and $t$, we conclude that $a_{s,h}$ is non-degenerate if and only if $h(s) \not \equiv 0 \pmod{r^2}$.

For the relationship with the hyperelliptic Ate pairing $a$, we use the fact that $t(D_2, D_1) = a(D_2, D_1)^{kq^{k-1}}$ \cite[Theorem 2]{GHV:2007}.

\end{proof}

\subsection{Examples of HV pairings}
\label{ss:HVexamples}

In this section, we describe how the pairings in the current literature  fit into the HV framework. While these pairings can be expressed as $a_{s,h}$ for some $s \in \Z$ and $h(x) \in \Z[x]$, their actual computation takes an alternate form in order to make use of shorter Miller loops.

\begin{enumerate}
\item
The {\it generalized Ate pairing} or {\it Ate$_i$ pairing},  was defined by Zhang \cite{Zhang:2008} as the analogue of the Ate$_i$ pairing for elliptic curves \cite{ZZH:2008-1}.
For $s \equiv q^j \pmod{r}$,
\begin{eqnarray*}
a_{s}: \G_2 \times \G_1 &\rightarrow& \mu_r \\
(D_2,D_1) &\mapsto& f_{s, D_2}(D_1)^{(q^k - 1)/r}.
\end{eqnarray*}
Since $r  \mid (q^k-1)$,  we may assume $0 < j < k$.  Note that if $s = q^j$ then no final exponentiation is needed, as is the case for the hyperelliptic Ate pairing. However, this choice of $s$ is never an improvement over
the Ate pairing as the Miller loop length is $i\log_2 q \geq \log_2 q$.

For $s \not\equiv q^j \pmod{r^2}$, it is straightforward to show this is the HV pairing $a_{s,h}$ where $h(x) = x-q^j$. Writing $s = q^j + \ell r$ for $\ell \in \Z$, we have
$(f_{s,D}) = sD - \rho(sD) = (f_{s,h,D}) + \ell rD$. As $\ell rD  \sim 0$, these functions differ only by a constant and thus give the same value after the final exponentiation.

\item The Ate pairings defined by Vercauteren \cite[Theorem 1]{Vercauteren:2008} for elliptic curves can be generalized directly to hyperelliptic curves.  To define the pairing, we
first choose an integer $m$ relatively prime to $r$ and express $mr$ in base $q$ as $mr = \sum_{i = 0}^n h_iq^i$. We can decompose the $m^{th}$ power of the Tate-Lichtenbaum pairing as
\begin{equation}
\label{ver}
t(D_2, D_1)^m = f_{\sum_{i = 0}^n h_iq^i, D_2}(D_1)^{(q^k-1)/r} = \left( \prod_{i=0}^{n} f_{h_iq^i, D_2}(D_1) \cdot \prod_{j = 0}^{n-1} g_j(D_1) \right)^{(q^k-1)/r} 
\end{equation}
where the $g_j$ ($j=0,\ldots,n-1$) are auxiliary functions defined through
\[
f_{\sum_{i = j}^n h_iq^i, D_2} = f_{\sum_{i=j+1}^{n} h_iq^i, D_2} f_{h_jq^j, D_2} g_j.
\]
The pairing $a_{[h_0,...,h_n]}$ is then defined as
\begin{eqnarray*}
a_{[h_0,...,h_n]}: \G_2 \times \G_1 &\rightarrow& \mu_r \\
(D_2,D_1) &\mapsto& \left( \prod_{i=0}^{n} f_{h_i, D_2}(D_1)^{q^i} \cdot \prod_{j = 0}^{n-1} g_j(D_1)  \right) ^{(q^k-1)/r}.
\end{eqnarray*}
It is easy to see
that  $a_{[h_0,...,h_n]}(D_2, D_1)$ equals  $t(D_2, D_1)^m$.
Indeed, by definition of the Miller functions and the action of the Frobenius on $D_1$ and $D_2$, we have that 
$$f_{h_iq^i, D_2}(D_1) = f_{q^i, D_2}(D_1)^{h_i} f_{h_i, q^iD_2} (D_1) = f_{q^i, D_2}(D_1)^{h_i} f_{h_i, D_2} (D_1)^{q^i}$$
as in the proof of \cite[Theorem 1]{Vercauteren:2008}.  While not explicitly noted in that proof, it is also true that $$ \left( \prod f_{q^i, D_2}(D_1)^{h_i} \right)^{(q^k-1)/r}= 1,$$
by an argument similar to that of Theorem \ref{HVtheorem}.
Therefore $a_{[h_0,...,h_n]}(D_2, D_1) = t(D_2, D_1)^m$. 
Thus, this pairing is simply the HV pairing $a_{q, h}$ where $h(x) = \sum_{i = 0}^n  h_i x^i$. 

The pairing $a_{[h_0,...,h_n]}$ is computed as a product of many Miller functions, as well as the auxiliary functions, and the total sum of the lengths of the Miller loops of the functions is $\sum_{i = 0}^n \log_2 h_i$. Thus, for efficiency, this pairing is fastest if the coefficients of $mr$ in base $q$ expansion are small. Vercauteren gives an algorithm to find suitable multiples of $r$ by searching for shortest vectors in a lattice spanned by vectors involving powers of $q$ \cite[\S3.3]{Vercauteren:2008}. This is the ``lattice" idea which was further generalized by Hess \cite{Hess:2008}. See Section \ref{ss:optimalpairings} for a discussion of the smallest loop length possible.

\item
The {\it R-ate pairing}, introduced by Lee, Lee and Park in 2008 \cite{LLP:2009}, was the first pairing defined as a ratio of generalized Ate pairings. 
We give a specific instantiation as an example (cf.\ \cite[Corollary 3.3(3)]{LLP:2009}). Let $T_i \equiv q^i\pmod r$ and $T_j\equiv q^j\pmod r$, where $0<i<j<k$, and write $T_i = aT_j + b$ for some $a,b \in \Z$. Then the R-ate pairing is
\begin{eqnarray*}
R: \G_2 \times \G_1 &\rightarrow& \mu_r \\
(D_2,D_1) &\mapsto& \left( f_{a,T_j D_2}(D_1) f_{b,D_2}(D_1) g(D_1) \right)^M, 
\end{eqnarray*}
where $g$ is an auxiliary function with divisor $aT_jD_2 + bD_2 - \rho(aT_jD_2 + bD_2)$ and $M \in \mathbb{N}$ is a final exponent. (The function $g$ is the analogue of the ratio of a linear and vertical function for the elliptic curve case.) 
Although it is ambigous in the original paper, this pairing requires a final exponentiation to yield a unique value. The exponent $M =(q^k-1)/r$ is sufficient, though a smaller exponent may also work, depending on the multiplicative orders of $T_i$ and $T_j$ modulo $r$ (see \cite[Corollary 3.3(3)]{LLP:2009} for details).
It is easy to work out (\cite[Theorem 3.2]{LLP:2009}) that 
\[
R(D_2, D_1) = \left( f_{T_i,D_2}(D_1)/f_{T_j,D_2}(D_1)^a \right)^M, 
\]
and thus $R$ is in fact a ratio of generalized Ate pairings.
Since $f_{a,T_j D_2}(D_1) = f_{a,D_2}(D_1)^{q^j}$ (\cite[Theorem 1]{ZZH:2008-1}),  in practice, the R-ate pairing is  computed as
\[
R(D_2, D_1) = \left( f_{a,D_2}(D_1)^{q^j} f_{b,D_2}(D_1) g(D_1) \right)^M.
\]
In this form, and with $M = (q^k-1)/r$, it is a straightforward calculation to establish that $R$ corresponds to the above Vercauteren pairing $a_{[h_0, ..h_i,., h_j]}$ with $h_0 = b, h_i = -1$, $h_j = a$ and all other coefficients equal to zero: let $\ell_i, \ell_j\in\Z$ such that $T_i = q^i + \ell_ir$ and $T_j = q^j + \ell_jr$, and express the $r$-multiple $(\ell_i - a\ell_j)r$ in base $q$, and use that $f_{1, D_2}$ is a constant function and therefore eliminated by the final exponentiation.
 In other words, $R$ is the HV-pairing $a_{s,h}$ where $s = q$ and $h(x) = ax^j - x^i + b$. 
\end{enumerate}

\subsection{Twisted Ate pairing}
\label{ss:twistedate}

In this section, we discuss the {\it twisted Ate pairing} $e: \G_1 \times \G_2 \rightarrow \mu_r$. The twisted Ate pairings use the fact that in certain situations, there is a ``twist" of the Frobenius $\pi$ which acts as $[q]$ on $\G_1$ and $[1]$ on $\G_2$, thereby reversing the roles of these groups in the Ate pairing. The main benefit of such pairings is that $D_1 \in \G_1$ is defined over $\F_q$, which means computing the Miller function $f_{s,D_1}$ is simpler. An added benefit is that the points in $D_2 \in \G_2$ have $x$-coordinates in a subfield of $\F_{q^k}$ which also may simplify the evaluation, as explained in Section \ref{ss:FinalExponentiation}.

Let $\C$ be a curve over a finite field $K = \F_q$. A \emph{twist} of $\C$ is a curve $\C'$ over $\F_q$ such that there exists an isomorphism $\phi: \C' \rightarrow \C$ defined over $\F_{q^{\delta}}$ for some ${\delta} \in \Z^+$.
If ${\delta}$ is the minimal degree extension of $\F_q$ over which the isomorphism is defined, then the twist $\C'$ is of {\it degree ${\delta}$.}  For more on twists of curves, see Silverman \cite[\S10.2]{Silverman:1992}.

Let $\pi$ be the Frobenius of $\C$ and let $\phi^\pi$ denote the isomorphism $\C' \rightarrow C$ obtained by $\pi$ acting on the coefficients of $\phi$. Then $\phi^\pi \circ \phi^{-1}$ is an automorphism of $\C$ of order $\delta$ in $\Aut(\C)$. Thus to look at twists of $\C$, one needs to consider the automorphism group of $\C$.  For genus 2 hyperelliptic curves over $\F_q$, $\Aut(\C)$ is isomorphic to one of the following groups  \cite{Cardona:2003,Cardona:2005}:
\[
\label{allauts}
C_2, C_{10}, C_2 \times S_3, V_4, D_8, D_{12}, 2D _{12}, \tilde{S_4}, \tilde{S_5}, M_{32}, \,\text{or}\, M_{160},
\]
where $C_n$ is the cyclic group of order $n$, $V_4$ is the Klein 4-group, $D_n$ is the dihedral group of order $n$, $S_n$ is the symmetric group of order $n$, $M_n$ is the group of order $n$ arising from a certain exact sequence \cite[Equation 6]{Cardona:2005}, and $2D_{12}, \tilde{S_4}, \tilde{S_5}$ are 2-coverings of $D_{12}$, $S_4$, and $S_5$, respectively. This implies that $\delta$, as the order of an element in $\Aut(C)$, has to divide $\#\Aut(C)$ for one of the above automorphism groups.

If $\C$ has a twist of degree $\delta$ with $m = \gcd(k, \delta) > 1$, then it is possible to define a non-degenerate, bilinear pairing on $\G_1 \times \G_2$. For applications to cryptography, we are interested in using the highest degree twist available, because elements of $\G_2$ can then be represented as elements of the Jacobian of the twist $\C'$ defined over $\F_{q^{k/m}}$. 

Given a curve $\C$, let $r \mid \#\Jac_C(\F_q)$ be a large prime, $k$ the embedding degree, and $C'$ a degree $\delta$ twist of $\C$.  We have an injection
\begin{align*}
[\cdot]: \mu_{\delta} &\ra \Aut(\C) \\
\xi &\mapsto [\xi]\end{align*}
where $\xi$ is the automorphism defined by the twist. Then $\G_2 = \jac_C(\F_q)[r] \cap \ker(\pi - [q]) = \jac_C(\F_q)[r] \cap \ker([\xi]\pi^{k/m} -1)$, and Zhang proved  the following theorem (\cite[Theorem 2]{Zhang:2008}):

\begin{theorem} 
Let $C$ be a hyperelliptic curve over $\F_q$ with a twist of degree $\delta$. Let $m = \gcd(k,\delta)$ and  $e=k/m$. 
Then
\begin{eqnarray*}
a^{\text{twist}}:\G_1 \times \G_2 &\rightarrow& \mu_r \\
({D_1},{D_2}) &\mapsto&  f_{q^e, D_1}(D_2),
\end{eqnarray*}

where  the representatives of $D_1 \in \G_1$ and $D_2 \in \G_2$ have disjoint support,  defines a non-degenerate bilinear pairing called the hyperelliptic twisted Ate pairing.
 \end{theorem}

\begin{remark}
For $\C$ with $\gcd(k,\#\Aut(\C)) \neq 1$, any pairing on $\G_1,\G_2 \subset \Jac_C(\F_q)$ in the HV framework has a twisted version, $a_{s,h}^{\text{twist}}: \G_1 \times \G_2 \rightarrow \mu_r$ \cite[Theorem 1]{Hess:2008}.
\end{remark}

We now define the {\it eta pairing}, which is essentially the twisted Ate pairing on supersingular curves, although historically it was introduced before the Ate pairing. The eta pairing makes use of  a
 {\it distortion map} on $\C$ instead of a twist. Let $e(\cdot,\cdot)$ denote any bilinear, non-degenerate, Galois-invariant pairing on $\Jac_C(\F_q)[r].$ A non-degenerate pairing ensures that given a non-zero divisor class $D_1$ of order $r$, there exists $D_2$ such that $e(D_1,D_2) \neq 1$. However, there are certain instances where a specific $D_1$ and $D_2$ pair to 1, for example, where $D_1, D_2$ both are defined over $\F_q$ and the embedding degree $k > 1$. To remedy this, we introduce distortion maps.
\begin{definition}Let $e$ be a non-degenerate pairing and $D_1$ and $D_2$ non-zero divisor classes of prime order $r$ on $\C$. A \emph{distortion map} is an endomorphism $\psi$ of $\Jac_C(\F_q)$ such that $e(D_1,\psi(D_2)) \neq 1$.\end{definition}

Galbraith et al. \cite{GPRS:2009} proved that distortion maps always exist for supersingular abelian varieties:
\begin{theorem}Let $A$ be a supersingular abelian variety of dimension $g$ over $\F_q$, and let $r$ be a prime not equal to the characteristic of $\F_q$. For every two non-trivial elements $D_1$ and $D_2$ of $A(\F_q)[r]$, there exists an endomorphism $\psi$ of $A$ such that $e(D_1,\psi(D_2))\neq 1.$
\end{theorem}

The {\it eta pairing} has been introduced in 2007 by Barreto et al. \cite{BGOS:2007} for supersingular curves.
It provides a generalization of the results of Duursma and Lee \cite{DuursmaLee:2003} for a specific instance of supersingular curves. Consider a supersingular curve $\C/\mathbb{F}_q$ (having one point at infinity) which has even embedding degree $k>1.$ Let $D_1$ and $D_2$ be reduced divisors of degree zero on $\C$ defined over $\mathbb{F}_q$ representing divisor classes with order $r.$ Assume that there exists a distortion map $\psi$ which allows for {\it denominator elimination} (see Section \ref{ss:FinalExponentiation}), meaning the $x$-coordinates of points in $\psi(D_2)$ lie in a subfield of $\F_{q^k}$. \begin{definition}
For $T\in\mathbb{Z},$ the \emph{eta pairing} $\eta_T$ is given by
\begin{eqnarray*}
\eta_T:\G_1\times\G_1 &\rightarrow& \mu_r \\
(D_1,D_2) &\mapsto&  f_{T,D_1}(\psi(D_2))^{(q^k - 1)/r}.
\end{eqnarray*}
\end{definition}
Note that in the literature, the eta pairing is often defined without 
the final exponent, though it is necessary to obtain a unique value in $\mu_r$.
In general, this pairing is not a non-degenerate, bilinear pairing, but Barreto et al. \cite[Theorem 1]{BGOS:2007} give sufficient conditions on $T$ under which $\eta_T(\cdot,\cdot)$ can be related to the modified Tate-Lichtenbaum pairing. In particular, this implies that for certain values of $T$, the eta pairing is indeed non-degenerate and bilinear. Moreover, the recent work of Lee, Lee and Lee \cite{LeeLee:2008, LeeLeeLee:2008} allows us to compute the eta pairing on genus $2$ curves for general divisors, which lifts an earlier restriction to the case of degenerate divisors  (see Section \ref{ss:degeneratedivisors}).

%*******************************************************

\section{Fast Computation of Hyperelliptic Pairings}
\label{s:fast}

In this section, we summarize the state of the art for fast computation of pairings on hyperelliptic curves of genus 2.

\subsection{Miller's algorithm}
\label{subsec: miller}

The algorithm used to compute Weil and Tate-Lichtenbaum pairings on elliptic curves was devised by Victor Miller in 1985 \cite{Miller:2004} and can be adapted to all pairings discussed in this paper \cite{ELM:2004}. Referring to the pairing definitions of Section \ref{s:pairings} one sees that to compute a pairing, it is necessary to evaluate a Miller function at a divisor. Algorithm \ref{algorithm: miller}, futheron referred to as ``Miller's algorithm", computes such a value using the structure of an addition chain for $s$.

Usually, an addition chain takes the form of a double-and-add chain, as follows.  Starting with the integer $k=0$, at each step one performs one of two possible calculations to update the value of $k$:  one either doubles to obtain $k \rightarrow 2k$ or doubles-and-adds to obtain $k \rightarrow 2k+1$.  To determine the sequence of steps needed to obtain any desired integer $s$ in this way, one reads the binary digits of $s$ from left to right, doubling once for each `0' and doubling-and-adding for each `1.' (For example, $5 = 101_2$ is obtained as $0 \rightarrow 2(0)+1 = 1 \rightarrow 2(1) = 2 \rightarrow 2(2)+1 = 5$.)  Starting from $0$, this algorithm computes $s$ in $\lfloor \log_2 s \rfloor +1$ steps (each of which consists of either one or two additions).

Miller's Algorithm computes $f_{s,D}$ following this double-and-add process by computing the Miller function $f_{k,D}$ at each step along the way, obtaining $f_{s,D}$ at the end. A double step involves one addition, and a double-and-add step involves two.  For each addition, we compute the new Miller function $f_{i+j,D}$ from the previously computed $f_{i,D}$ and $f_{j,D}$ via the relationship
\[
f_{i+j,D} = f_{i,D} f_{j,D} h_{iD,jD}, \qquad i,j>0,
\]
where the auxiliary function $h_{D',D''}$ is a function with divisor
\[
\rho(D') + \rho(D'') - \rho(D' + D'').
\]
The computation of $h_{D',D''}$ is performed by an enhanced version of Cantor's Algorithm (cf. Section \ref{subsec: arithjac}), here Algorithm \ref{algorithm: millerstep}.  It is called under the name \texttt{Cantor()} once (if doubling) or twice (if doubling and adding) in each for-loop of Miller's Algorithm. 
Using the result of Algorithm 2, one calculates $f_{2i,D}$ from $f_{i,D}$ (``double'') or $f_{2i+1,D}$ from $f_{i,D}$ and $f_{1,D}$ (``double and add''),  where $f_{1,D}$ is a constant function.

In order to compute the pairing value, the Miller function $f_{s,D_2}$  must be evaluated a divisor $D_1$, but 
this evaluation is not possible unless $D_1$ and $D_2$ have disjoint support, which is not the case if both are reduced.
However, using reduced divisors and Mumford
representation is  too useful to dispense with, so the solution is
the following.  Let $z$ be a {uniformizer} at $P_\infty$ (for example, $z(x,y)= x^2/y$ is a convenient choice).
Then, if $f$ is a function with order $-r$ at $P_{\infty}$, define the {\em leading coefficient} at $P_\infty$ of $f$, denoted as $\lc_{\infty}(f)$, to be $(z^rf)(P_{\infty})$.
Then the \emph{normalization} of $f$ is the scalar multiple $f^{\,\norm} = f/\lc_\infty(f)$ which has leading coefficient $1$.  For the hyperelliptic Ate pairing \cite[Lemma 6]{GHOTV:2007}, when $z$ is $\mathbb{F}_q$-rational,
\[
f_{q,\rho(D_2)}(D_1) = f^{\,\norm}_{q,\rho(D_2)}(\epsilon(D_1)).
\]
The right-hand expression requires computing the leading coefficient, but solves the problem of non-disjoint supports of $D_1$ and $D_2$ without losing the usefulness of Mumford representation.

For HV Pairings and the modified Tate-Lichtenbaum pairing, the same solution is possible.  Consider the computation of $t(D_2, D_1) = f_{r, D_2}(D_1)^{(q^k-1)/r}$ where $D_1, D_2$ are reduced. 
Let $-b_i$ be the coefficient of $P_{\infty}$ in $D_i$ for $i = 1,2$. (Note that $b_i = -1$ or $-2$, depending on whether or not the reduced divisor $D_i$ is degenerate.)  
The function $f_{r, D_2}$ has divisor $rD_2$ with order $-b_2r$ at $P_\infty$. 
Therefore, if $z$ is an $\F_{q^k}$-rational uniformizer at $P_\infty$,  
\[
 f_{r, D_2}(D_1) =f^{\, \norm}_{r,D_2}(\epsilon(D_1)) / z(P_{\infty})^{b_1b_2r}.
\]
Since $b_1b_2r$ is a multiple of $r$, the contribution of $z(P_{\infty})^{b_1b_2 r (q^k-1)/r}$ is 1, and thus 
\[
 f_{r, D_2}(D_1)^{(q^k-1)/r}=f^{\, \norm}_{r,D_2}(\epsilon(D_1))^{(q^k-1)/r}. 
\]
As the HV pairing $a_{s,h}(D_2, D_1)$ is a simply a power of the modified Tate pairing $t(D_2, D_1)$ (see Theorem \ref{HVtheorem}),  in whichever form the pairing $a_{s,h}(D_2, D_1)$ is computed, evaluating normalized functions at effective divisors will give the pairing value. 

In the elliptic curve case,  it is more efficient to evaluate the Miller functions and the auxiliary functions $h_{D',D''}$ at the desired divisor (denoted $D_2$ in Miller's Algorithm) at each step, instead of reserving the evaluation for the end.  In order to allow for this, $D_2$ is passed to Cantor's Algorithm. 
We now turn to a discussion of this aspect in the case of hyperelliptic curves.

In Miller's Algorithm, the current Miller function $f$ is stored as two polynomials $f_1$ and $f_2$ such that $f = f_1/f_2$.  Similarly, the auxiliary functions $h$ are returned from Cantor's Algorithm as $h_1$ and $h_2$.  It remains to explain how to evaluate a polynomial function $g(x,y)$ on $\C$ at the effective part of a divisor given in Mumford representation $(u(x), v(x))$ (we need only the effective part because of the preceeding discussion and the computation of the leading coefficient).  We need to evaluate $G(x) = g(x, v(x))$ at the zeroes of $u(x)$.  This is the same as computing the resultant $\operatorname{Res}(G(x), u(x))$.
 Performing a resultant calculation is sufficiently costly that it is best left to the end of Miller's Algorithm, as long as the size of the Miller functions can be kept low in the meantime.  Fortunately, in preparation for the eventual final resultant, it suffices to compute the Miller functions in $x$ and $y$ modulo $u(x)$, while substituting $y=v(x)$, effectively capping their degrees.

\begin{algorithm}
\caption{Miller's Algorithm}
\label{algorithm: miller}
\begin{algorithmic}[1]
\REQUIRE $D_1 = (u_1,v_1)$, $D_2 = (u_2,v_2)$, $d$, $s = \sum_{i=0}^N s_i 2^i$ \\
\ENSURE $f^{\,\norm}_{s,D_1}(\epsilon(D_2))^d$ \\
\STATE $D \leftarrow D_1$ \\
\STATE $f_1 \leftarrow 1, f_2 \leftarrow 1, f_3 \leftarrow 1$ \\
\FOR{$i= N-1$ down to $0$}
 \STATE $f_1 \leftarrow f_1^2 \pmod{u_2}, f_2 \leftarrow f_2^2 \pmod{u_2}, f_3 \leftarrow f_3^2$ \\
 \STATE $(D, h_1, h_2, h_3) \leftarrow \operatorname{Cantor}(D,D,D_2)$ \\
 \STATE $f_1 \leftarrow f_1 \cdot h_1 \pmod{u_2}, f_2 \leftarrow f_2
\cdot h_2 \pmod{u_2}, f_3 \leftarrow f_3 \cdot h_3$ \\
 \IF{$s_i = 1$}
   \STATE $(D, h_1, h_2, h_3) \leftarrow \operatorname{Cantor}(D,D_1, D_2)$ \\
   \STATE $f_1 \leftarrow f_1 \cdot h_1 \pmod{u_2}, f_2 \leftarrow f_2
\cdot h_2 \pmod{u_2}, f_3 \leftarrow f_3 \cdot h_3$ \\
 \ENDIF\\
\ENDFOR\\
\STATE {\bf return} $\left( \operatorname{Res}(f_1,u_2)/(f_3^{\deg(u_2)}\cdot \operatorname{Res}(f_2,u_2)) \right)^d$ \\
\end{algorithmic}
\end{algorithm}

\begin{algorithm}
\caption{Cantor's Algorithm}
\label{algorithm: millerstep}
\begin{algorithmic}[1]
\REQUIRE  $D_1 = (u_1,v_1)$, $D_2 = (u_2,v_2)$, $D' = (u,v)$ \\
\ENSURE $\rho(D_1+D_2), f(x,v(x)) \pmod{u}, g(x,v(x)) \pmod{u}, \lc_\infty(h_{D_1, D_2})$ where $h_{D_1,D_2} = f/g$ \\
\STATE compute $(d_1,e_1,e_2)$ such that $d_1 = e_1u_1+e_2u_2 =
\gcd(u_1,u_2)$ \\
\STATE compute $(d,c_1,c_2)$ such that $d = c_1d_1+c_2(v_1+v_2+H) =
\gcd(d_1,v_1+v_2+H)$ \\
\STATE $s_1 \leftarrow c_1e_1, s_2 \leftarrow c_1e_2, s_3 \leftarrow c_2$ \\
\STATE $U \leftarrow (u_1u_2)/d^2, V \leftarrow (s_1u_1v_2+s_2u_2v_1 +
s_3(v_1v_2+F))/d \pmod{U}$ \\
\STATE $f \leftarrow d \pmod{u}, g \leftarrow 1, h \leftarrow 1$ \\
\WHILE{$\deg(U) > g$}
 \STATE $U' \leftarrow (F-VH - V^2)/U, V' \leftarrow (-H-V) \pmod{U'}$ \\
 \STATE $f \leftarrow f \cdot (v - V) \pmod{u}$ \\
 \STATE $g \leftarrow g \cdot U' \pmod{u}$ \\
 \IF{$\deg(V) > g$}
   \STATE $h \leftarrow -\operatorname{leadingcoeff}(V) \cdot h$\\
 \ENDIF\\
 \STATE $U \leftarrow U', V \leftarrow V'$ \\
\ENDWHILE\\
\STATE {\bf return} $(U,V), f,g,h$\\
\end{algorithmic}
\end{algorithm}

If Steps 5 and 8 through 13 are removed from Cantor's Algorithm and only $(U,V)$ is returned, the algorithm computes
$\rho(D_1+D_2)$ for any divisors $D_1$ and $D_2$ in Mumford
representation (this is the usual meaning of ``Cantor's Algorithm'' as in Section \ref{subsec: arithjac}).  If these steps are included, then Cantor's Algorithm
 can also return $f, g \pmod{u}$ such that
$f/g = h_{D_1,D_2}(x,v(x))$ for some specified divisor $(u,v)$.  This is the form in which it is used in Miller's Algorithm.

In the case that we are pairing degenerate divisors (see Section
\ref{ss:degeneratedivisors}), a norm computation may be
preferred to the resultant method \cite{GHV:2007}.

\subsection{Using effective divisors and the leading coefficient}

The leading coefficient of $f_{s,D}$ is an element of the field of definition of the function. Therefore, in the case of twisted pairings, the leading coefficient of $f_{s, D_1}$ is defined over $\F_q$. Therefore,
if the pairing includes a final exponentiation, the leading coefficient will be eliminated and thus may be ignored in the computation of the pairing.

\subsection{Final exponentiation}
\label{ss:FinalExponentiation}

As described in Section \ref{s:pairings}, most of the hyperelliptic
pairings involve a {\it final exponentiation} of a Miller function
$f_{s,D}(D')$ by $(q^k-1)/r$,
where $D \in \Jac_C(\F_q)[r]$ and $D'$ is an arbitrary divisor in
$\Jac_C(\F_{q^k})$.
As has been widely reported, this extra computation
has its benefits, in particular when $k$ is even. Many of these are described by Scott
\cite{Scott:2007} and Galbraith, Hess, and Vercauteren \cite{GHV:2007}; we summarize the main
ones here.

When $k$ is even, the field $\F_{q^k}$ can be constructed as a degree two extension of $\F_{q^\ell}$, where $2\ell = k$.
We can represent elements  as
$a + ib$ with $a,b \in \F_{q^\ell}$ and $\gamma^2$  a quadratic non-residue over
$\F_{q^\ell}$. It is straightforward to check that
\[( 1/(a + \gamma b))^{q^\ell - 1}  = (a - \gamma b)^{q^\ell - 1}\]
which means  inversion can be replaced by conjugation since the result is the same after final
exponentiation. In particular, this applies to any denominators of computations in  Miller's algorithm.

There is a further optimization, {\it denominator elimination}, which in fact allows one to ignore all denominators in Miller's algorithm. In computing $f_{s,D}(D')$ where $D$ is a divisor defined over the base field $\F_q$,
one computes the numerator and denominator values separately (see Algorithm \ref{algorithm: miller}). If $D' = (u(x), v(x))$ has $u(x)$ defined over $\F_{q^\ell}$,
then the computation of the denominator involves only $D$ and $u(x)$ and therefore  becomes trivial after final exponentiation.
In the case of supersingular curves, for example, a suitable evaluation divisor can be found using a distortion map $\psi$ (see Section \ref{ss:twistedate}) such that
$\psi(D')$ has $x$-coordinates in $\F_{q^\ell}$ \cite{GPRS:2009}.

The final exponentiation is generally computed in multiple steps by
writing $(q^k - 1)/r$ as a product of  polynomials in base $q$ expansion and
exploiting finite field constructions, in particular the $q^{\th}$ power of
Frobenius, which speeds up computation \cite{GHV:2007}. 
Other methods for faster computation include signed sliding window methods
\cite{GPS:2006a}, as well as trace and tori methods \cite{GalbraithScott:2008,GPS:2006b}.

\begin{remark}
\label{rmk:ate} As the Ate pairing does not require final exponentiation, these techniques are unavailable. Furthermore,
as stated by Granger et al., there are also possible security implications; namely, the problem of {\it pairing inversion} (given $\gamma$ and $D_1$, find $D_2$ such that $a(D_1,D_2) = \gamma$) may not be as hard (see \cite[Intro.]{GHOTV:2007}). However, we remark that if $r^2 \nmid (q^k-1)$ and $r$ is prime, a superfluous final exponentiation of the Ate pairing still gives a non-degenerate result.
\end{remark}

\subsection{Degenerate divisors}
\label{ss:degeneratedivisors}
For a genus 2 curve, a general reduced divisor $D$ is of the form $D = (P_1) + (P_2) -2(\infty)$ and a degenerate divisor is of the form $D = (P) - (\infty)$. As there are fewer points in the support, the arithmetic is faster  when adding a general divisor to a degenerate divisor than when adding two general
divisors. This speeds up the computation of the Miller function $f_{s,D}$ where $D$ is degenerate. Furthermore, the evaluation of a Miller function on a degenerate divisor is also faster by at least half, since there is only one affine point. Many of the fastest hyperelliptic pairing computations use degenerate divisors, including the examples noted with
[$a$], [$b$] and [$c$] in the Table \ref{pairingtable}. We summarize here when it is possible to use degenerate divisors as either the first or second argument of a pairing.

Should $\Jac_C(\F_q)$ be of prime order $r$, then for any $P \in C(\F_q)$, the divisor $D = (P) - (\infty)$ can be used as the first argument, regardless of the pairing.
Furthermore, if $C$ is supersingular, then using a distortion map $\psi$ (see Section \ref{ss:twistedate}), we have that $\psi(D)$ is also degenerate and pairs non-trivially with $D$. Hence, for supersingular curves with prime-order $\Jac_C(\F_q)$, we can use degenerate divisors as both arguments of the Tate-Lichtenbaum pairing. This fact was originally exploited in the definition of the $\eta_T$ pairing by Duursma and Lee \cite{DuursmaLee:2003}.
In the more general situation where $\#\Jac_C(\F_q)$ is not prime and/or the curve $C$ is not supersingular, using degenerate divisors is not as straightforward, as noted by Frey and Lange \cite{FreyLange:2006}.
If $\#\Jac_C(\F_q) = nr$ where $\gcd(n,r) = 1$, there is no guarantee that there exists a degenerate divisor $D$ of order $r$. The probability that a reduced divisor is of order $r$ is $1/n$ and the probability  that a divisor is degenerate is roughly  $1/q$, by the Hasse-Weil bounds on $C(\F_{q})$ and $\Jac_C(\F_q)$. Therefore, assuming independence, a heuristic argument gives that
the probability a divisor is degenerate and order $r$ is $1/qn$. This implies that using a degenerate divisor for the first argument is not necessarily possible.

However, Frey and Lange \cite{FreyLange:2006} show that for $q$ large enough (as in a cryptographic setting), it is possible to use a degenerate divisor as the second argument. In other words, there exists $D_2 = (P) - (\infty) \in \Jac_C(\F_{q^k})$ such that for any $D_1 \in \Jac_C(\F_q)[r]$, the Tate-Lichtenbaum pairing $\tau (D_1, D_2)$ is non-trivial. 
The probability that $P \in C(\F_{q^k})$ yields such a divisor $D_2$ has a lower bound of $1/k\log_2 q.$ Moreover, if $k = 2d$ is even, it is possible to choose $P = (x,y)$ with $x \in \F_{q^d}$ and $y \in \F_{q^k}$, using a degenerate divisor on the quadratic twist of $C/\F_{q^d}$.
This technique is used for example by Fan, Gong and Jao \cite{FGJ:2008} and allows for denominator elimination.

\begin{remark}
As remarked by Galbraith, Hess and Vercauteren \cite[\S7]{GHV:2007}, there are potential security implications with using degenerate divisors, depending on the application.  While the discrete logarithm problem with a degenerate divisor as a base point is no easier than that with a general divisor \cite{KKAT:2004}, other hardness assumptions such as pairing inversion (see Remark  \ref{rmk:ate}) are potentially compromised, as Granger et al. have noted \cite{GHOTV:2007}. To our knowledge, the topic remains unresolved.

We also remark that there are protocols in which it may not always be possible to use degenerate divisors, for example, when computing a pairing where one input is required to be a random multiple of a divisor $D$.
\end{remark}

\subsection{Rubin-Silverberg point compression}\label{ss:rubinsilverberg}

Another method available to us in genus 2 is the point compression technique of Rubin and Silverberg \cite{RubinSilverberg:2002}, who note
that supersingular abelian varieties can be identified with
subvarieties of Weil restrictions of supersingular elliptic
curves.

Recall that a \emph{supersingular $q$-Weil number} is a complex
number of the form $\sqrt{q}\zeta$, where $\zeta$ is a root of
unity and $\sqrt{q}$ denotes the positive square root. Let $m$ be
the order of $\zeta$.

The following theorem allows us to define a useful invariant:
\begin{theorem}[\cite{RubinSilverberg:2002}]
Suppose $A$ is a simple supersingular abelian variety of dimension $g$ over $\F_q$, where $q$ is a power of a prime $p$, and $P(x)$ is the characteristic polynomial of the Frobenius endomorphism of $A$. Then  $P(x) = G(x)^e$, where $G(x) \in \Z[x]$ is a monic irreducible polynomial with $e=1$ or 2. All of the roots of $G$ are supersingular $q$-Weil numbers.
\end{theorem}
We call the roots of $G$ the \emph{$q$-Weil numbers} for $A$.

\begin{definition}
The \emph{cryptographic exponent} of $A$ is defined by

\begin{align*}
c_A = \left\{
\begin{aligned}
\frac{m}{2}\ \ \ \ \ \ \ \ \ \ , &\text{\ if $q$ is a square} \\
\frac{m}{\gcd(2,m)},  &\text{\ if $q$ is not a square}.
\end{aligned}
\right.
\end{align*}

Let $\alpha_A = c_A/g$; it is the \emph{security parameter} of $A$.
\end{definition}

Now let $\mathbb{F} \subset \mathbb{F'}$ be finite fields, $E$ an elliptic curve over $\mathbb{F}$, and let $Q \in E(\mathbb{F'})$. Recall that the trace from $\mathbb{F'}$ to $\mathbb{F}$ is given by $$\Tr_{\mathbb{F'}/\mathbb{F}}(Q) = \sum_{\sigma \in \Gal(\mathbb{F'}/\mathbb{F})}\sigma(Q).$$ Rubin and Silverberg prove the following result:
\begin{theorem}[\cite{RubinSilverberg:2002}]
Let $E$ be a supersingular elliptic curve over $\F_q$, $\pi$ a $q$-Weil number for $E$ ($\pi \not\in \Q$). Fix $r \in \N$ with $\gcd(r,2pc_E) = 1$. Then there is a simple supersingular abelian variety $A$ over $\F_q$ having the following properties.
\begin{enumerate}
\item $\dim A = \varphi(r)$.
\item For every primitive $r^{\th}$ root of unity $\zeta$, $\pi\zeta$ is a $q$-Weil number for $A$.
\item $c_A = rc_E$.
\item $\alpha_A = (r/\phi(r))\alpha_E$.
\item There is a natural identification of $A(\F_q)$ with the following subgroup of $E(\F_{q^r})$ : $$\{Q \in E(\F_{q^r}) : \Tr_{\F_{q^r}/\F_{q^{r/l}}}(Q) = 0\;\text{for every prime}\; l\mid r\}.$$
\end{enumerate}
\end{theorem}

This theorem can be thought of as a form of point compression for
supersingular elliptic curves. More concretely, the theorem allows us to replace the
Jacobian of a hyperelliptic curve $C$ over $\F$ with an elliptic curve
$E$ over an extension $\F'$ of $\F$, while still exploiting the per-bit
security gain of higher genus hyperelliptic curves. From a security standpoint, there
is no difference between working with $E(\F')$ and working with $\Jac_C(\F)$.
On the other hand, one needs fewer bits to represent divisors with
support in $C(\F)$ than to represent points in $E(\F')$.

As noted by Galbraith \cite{Galbraith:2001}, recent implementations \cite{BGOS:2007} indicate that pairings on elliptic curves with the Rubin-Silverberg compression are, in general, more efficient than using the pairings on Jacobians of hyperelliptic curves. However, it seems that Rubin and
Silverberg have initiated a promising investigation into the
arithmetic geometry of abelian varieties and its applications to
pairings. Much work remains to be done, in particular with respect to the
torsion structure of these varieties.

\subsection{A comparison of pairings}
\label{ss:looplength}

We conclude this section by summarizing in Table
\ref{pairingtable} all known variants of the Tate-Lichtenbaum pairing defined in Section \ref{s:pairings},
in terms of their loop length and whether or not there is a final exponent of $(q^k-1)/r$. Note that if there is a final exponent, in the case of even embedding degree $k$, this allows for the optimizations described in Section \ref{ss:FinalExponentiation}. The last column gives references to specific examples of curves of genus 2 in the literature for which the efficiency of the pairing has been analyzed, either theoretically, via implementation or both.

All pairings in Table \ref{pairingtable} except the Tate-Lichtenbaum pairing and the modified Tate-Lichtenbaum pairing are defined on $\G_2 \times \G_1$,
but if $\gcd(k, \#\Aut(\C)) \neq 1$, then there exist the twisted versions on $\G_1 \times \G_2$ which have the same final exponent and loop length.\\

\begin{minipage}[]{4in} 
\label{pairingtable}
\centerline{{\sc Table} \ref{pairingtable}. A comparison of pairings.}
\ \\[-3ex]
\centering

\begin{tabular}{|c|c|c|c|c|}

\hline
\multirow{2}{*}{} Pairing & Curves & Final   & Loop  & Examples \\
					& 		& Exponent & Length & for $g = 2$ \\
\hline\hline
Modified Tate & All &  Yes & $\log_2 r$ &  \cite{FGJ:2008}$^a$, \cite[\S5]{OhEigeartaighScott:2006}, \cite{ChoieLee:2004},     \\	
\hline
Ate \cite{GHOTV:2007} & All & No & $\log_2 q$  & \\
\hline
\multirow{2}{*}{Eta \cite{BGOS:2007}} & \multirow{2}{*}{Supersingular} & \multirow{2}{*}{Yes} &  Varies  &  \multirow{2}{*}{\cite{BGOS:2007}$^b$} \\
&&& $(\log_2 q)$ possible& \\

\hline
\multirow{2}{*}{HV \cite{Hess:2008, Vercauteren:2008}} & \multirow{2}{*}{All} & \multirow{2}{*}{Yes} & Varies  & \multirow{2}{*}{\cite[\S4]{Vercauteren:2008}$^c$}  \\
&&& $(\log_2 r)/\varphi(k)$ possible  & \\
\hline
\multirow{2}{*}{Ate$_i$ \cite{Zhang:2008}} & \multirow{2}{*}{All} &  \multirow{2}{*}{Yes} & $\log_2 (q^i \pmod{r})$ &  \multirow{2}{*}{\cite[\S5]{Zhang:2008}$^d$}  \\	
	&&& $(\log_2 r)/\varphi(k)$ possible & \\
	
\hline
R-ate \cite{LLP:2009} & All &  Yes & Varies &   \cite[\S5]{LLP:2009}$^e$, \cite[\S4,5]{GLM:2008}$^f$ \\
\hline
\end{tabular}
\end{minipage}

\begin{itemize}
\item[[$a$\!\!\!]] Fan, Gong and Jao use efficiently computable automorphisms to compute a power of the modified Tate-Lichtenbaum pairing on two Kawazoe-Takahashi families of
non-supersingular curves over prime fields.  This algorithm allows for a theoretical reduction  of up to one fourth in the length of the Miller loop ($\log_2 r$).
They implement this on curves over $\F_p$ where $p$ is a 329-bit prime and $k = 4$ and compare this with pairings on a supersingular curve defined over
$\F_p$ with $p$ a 256-bit prime and $k = 4$.
Using all known optimizations
 (degenerate divisors, encapsulated group operations, final exponentiation, fast field arithmetic), the pairing computation on the non-supersingular curve is about $55.8\%$ faster.

 \item[[$b$\!\!\!]]  This is one of the fastest known pairing implementations on a hyperelliptic curve and makes use of many optimizations including degenerate divisors and a special octupling formula.

 \item[[$c$\!\!\!]]   Vercauteren gives an example of a family of supersingular curves  with $k = 12$ such that the loop length is approximately $\log_2 r/\varphi(k)$.

 \item[[$d$\!\!\!]]   Zhang gives examples of Kawazoe-Takahashi curves  with $k = 8, 24$ such that the twisted Ate$_i$ pairing
has loop length approximately $\log_2 r/\varphi(k)$.

 \item[[$e$\!\!\!]] Lee, Lee and Park show that for supersingular curves the loop length can  theoretically be approximately $(\log_2 q)/2$. They also compute an example on a Duursma-Lee curve with $k = 5$, achieving a loop length 21$\%$ shorter than the Ate.

 \item[[$f$\!\!\!]] Galbraith, Lin and Mireles Morales \cite{GLM:2008} describe how to use the R-ate pairing on a real model of a hyperelliptic curve of genus 2 over $\F_p$ with $k = 6$. By using a distortion map $\psi$
on $\Jac_C(\F_p)[r]$ such that the image of $\G_1$ is in the $p$-eigenspace, $\G_2$, they are also able to make use of denominator elimination. They conclude that such pairings are theoretically
competitive with both pairings on certain elliptic curves with $k = 3$ and with hyperelliptic curves in the imaginary model with $k = 4$.

\end{itemize}

%*******************************************************

\section{Future Work on Hyperelliptic Pairings}
\label{s:faster}

In this section, we present possible areas for future work, expanding upon the list in the 2007 survey paper of Galbraith, Hess and Vercauteren \cite{GHV:2007}. 
 We list some newer problems, mention some recent advancements in the elliptic curve case which may find generalizations in pairings for $g \geq 2$, and conclude by revisiting the 2007 list \cite{GHV:2007}.

\subsection{Achieving optimal loop length}
\label{ss:optimalpairings}
Since 2007, there has been a flurry of new work to reduce the loop length in Miller's algorithm using variants of the Ate pairing. In particular, the Ate pairing on hyperelliptic curves of genus $g$
already reduces the loop length by up to a factor of $g$ when compared to the Tate-Lichtenbaum pairing \cite{GHOTV:2007}. Vercauteren \cite{Vercauteren:2008}
uses the following definition to characterize pairings with certain loop lengths.
\begin{definition}  \cite{Vercauteren:2008} Let $e: \G_1 \times \G_2 \mapsto \mu_r \subset \F_{q^k}^*$ be a non-degenerate, bilinear pairing defined using a
combination of Miller functions. We call $e( \cdot, \cdot)$ an {\it optimal} pairing if it can be computed using  $(\log_2 r )/ \varphi(k) + \varepsilon(k)$ Miller iterations,
where $\varphi$ is the Euler phi function and $\varepsilon(k) \leq \log_2 k$.
\end{definition}
Note that this means a pairing is optimal if the total sum of all the loop lengths of the Miller functions is approximately $(\log_2 r )/ \varphi(k)$. 

For an HV pairing $a_{s,h(x)}$ with $h(x) = \sum_{i=0}^n h_ix^i$, the total sum of loop lengths is $\sum_{i=0}^n \log_2 h_i$. Thus to be optimal, 
it is necessary but not sufficient that the coefficients of $h$ are bounded by $r^{\varphi(k)}$. This can be achieved by finding the 
shortest vectors in a lattice spanned by vectors involving powers of $s$ \cite[\S3.3]{Vercauteren:2008}. Vercauteren and Zhang 
both give examples of genus 2 HV pairings (see Table \ref{pairingtable}) where the polynomial $h(x)$ satisfies this bound and has only one coefficient which is not $\pm 1$, 
therefore providing examples of optimal hyperelliptic pairings. It remains open whether given a hyperelliptic curve it is always possible to construct an optimal HV pairing.
One direction would be to look at extending the method of Vercauteren   \cite{Vercauteren:2008}  which constructs optimal pairings on {\it parameterized} families of  elliptic curves. 

Vercauteren also conjectures that for elliptic curves without  efficiently computable automorphisms other than the Frobenius, no pairing can be better than optimal \cite[\S 2]{Vercauteren:2008}. More specifically, he conjectures that for  such a curve, any non-degenerate pairing  requires at least $(1 - \delta)\log_2 r /\varphi(k)$ Miller iterations where $0 < \delta < 1/4$. For a curve with a set of efficiently computable endomorphisms $\mathcal{E} \subset \rm{End}(E)$, Vercauteren defines a {\it superoptimal} pairing as one which can be computed using 
$(\log_2 r )/ \# \mathcal{E} + \varepsilon(k)$ Miller iterations.  It remains to examine what is the best possible for genus 2 curves, both with and without the existence of efficiently computable endomorphisms (see also Section \ref{ss:automorphisms_faster}). Furthermore, it is not known whether there are other non-degenerate, bilinear hyperelliptic pairings on $\G_1 \times \G_2$ which are not part of the HV framework. 

Lastly, we remark that the computation of an HV pairing cannot be measured
solely by the sum of loop lengths.
There  is also the cost of computing the auxiliary functions (see (2),(3)
in Section \ref{ss:HVexamples}). It remains to formally compare the
cost of  these additional computations with the benefit of a shorter
total sum of  Miller loop lengths.

\subsection{Using efficiently computable automorphisms}
\label{ss:automorphisms_faster}

One newer method to speed up computations is to use efficiently computable automorphisms of the curve $\C$ (beyond the Frobenius). For example, 
Fan, Gong and Jao use efficiently computable automorphisms in computing a power of the modified Tate-Lichtenbaum pairing on some specific non-supersingular genus 2 curves over prime fields \cite{FGJ:2008}.
An open task is to explore how far can this be generalized to other genus 2 curves.

Furthermore, Hess \cite{Hess:2008} extends his pairing framework for ordinary elliptic curves 
to exploit efficiently computable automorphisms. This does not generally give an improved loop length since $\# \mathcal{E} \leq \varphi(k)$ for most ordinary elliptic curves.  However, as hyperelliptic curves have a greater variety of $\Aut(\C)$, it would be worthwhile to  examine what
improvements in loop length can be made by extending the HV framework to exploit these automorphisms.

\begin{comment}
In the case of ordinary elliptic curves, Hess \cite{Hess:2008} uses automorphisms of the elliptic curves to define {\it extended pairings} which are  variations of the HV pairings of Theorem \ref{HVtheorem}.  
Specifically, if $n$ is an integer dividing the least common multiple of $k$ and $\# \Aut(E)$, then it is possible to find $s \in \Z$ and $h(x)$ such that $a_{s,h}$ is non-degenerate and $\left\|h\right\|_1 = 
O(r^{1/\varphi(n)})$. Thus if  $\varphi(n) > \varphi(k)$, this yields a pairing function with a smaller degree. For ordinary elliptic curves, there are only a few cases where $\varphi(n) > \varphi(k)$, so this is not widely applicable. However, as hyperelliptic curves have a greater variety of $\Aut(\C)$, it would be interesting to examine what improvements 
could be made in those cases. 
\end{comment}

\subsection{Fast arithmetic and the embedding degree}
\label{ss:odddegree}

In the case of even embedding degree $k$, it is traditional to exploit
the degree two subfield, as explained in Section
\ref{ss:FinalExponentiation}. In fact, Koblitz and Menezes define
\emph{pairing friendly fields} to be finite fields of the form
$\mathbb{F}_{q^k}$ such that $k=2^i3^j$ for $0 \leq i,j \in \mathbb{Z}$
and $q \equiv 1 \pmod{12}$ \cite[\S5]{KoblitzMenezes:2005}. (If $k$ is
strictly a power of 2 then it is only required that $q \equiv 1 \pmod 4.$)
By a theorem of Lidl and
Niederreiter \cite[Theoreom 3.75]{LidlNiederreiter:1997} and more particularly, by
a specific instance of this theorem given by Koblitz and Menezes \cite[Theorem
2]{KoblitzMenezes:2005}, we can construct the extension $\mathbb{F}_{q^k}$
for $k$ of this form using a tower of quadratic and cubic extensions.
There are thus certain advantages we can make use of for $k=2^i3^j.$  For
instance, there exist fast arithmetic methods for degree 2 and 3
subextensions; namely, the Karatsuba method for quadratic subextensions and
the Toom-Cook method for cubic subextensions \cite[\S4.3.3]{Knuth:1997}.
These methods are used to economize the arithmetic in the smaller
fields which reduce the number of field multiplications. However, there are embedding degrees not of this form, particularly among recent constructions of non-supersingular curves,  
and hence it would be worthwhile to see if these ideas can be extended to embedding degrees $k$ containing other prime factors.

\subsection{Degenerate divisors} 
\label{ss:degeneratedivisors_faster}

As discussed in Section \ref{ss:degeneratedivisors}, one common optimization is to use degenerate divisors. Frey and Lange \cite{FreyLange:2006} give a lower bound on the probability that $P \in C(\F_{q^k})$ gives a non-trivial pairing value when used as a degenerate divisor in the second argument of the Tate-Lichtenbaum pairing.  However, to our knowledge, there is no method to efficiently find such points beyond 
simple trial and error.

We also consider using degenerate divisors with Ate-type pairings $a$ on $\G_2 \times \G_1$ (or twisted Ate on $\G_1 \times \G_2$).
While a heuristic argument shows that the likelihood that a divisor of $\G_1$ is degenerate is small, it would be useful to know if there are particular curves where this is more likely and if so, how to find such divisors. It also remains to analyze the likelihood that an element of $\G_2$ is degenerate. We note that for $D \in \G_2$, if $D = (P) - (\infty)$, then $\pi(D) = ( \pi(P) ) - (\infty)$ implies that the divisor class $qD$ is also degenerate.

\subsection{Ignoring the last bit} 
\label{ss:lastbit}

In the case of the modified Tate-Lichtenbaum pairing on elliptic curves, when computing $f_{r,D_1}(D_2)$, it is possible to ignore the last bit in the expansion of $r$. 
This follows from the fact that since $r$ is odd, the last iteration of the Miller loop of the Tate-Lichtenbaum pairing is the evaluation at $D_2$  of the line function corresponding to the line through $(r-1)P$ and $P$. 
This is a vertical line and so by the choice of divisor $D_2$ with $x$-coordinates lying over $\F_{q^d}$, this is eliminated by the final exponentiation. 
While this does not give a large improvement compared to other loop length reductions,  it is worth verifying whether this trick might be used in the case of hyperelliptic curves. 

\subsection{Compression and higher degree twists}
\label{ss:compression_faster}
Galbraith and Lin \cite{GalbraithLin:2009} give explicit formula to compute the Weil pairing on elliptic curves given only $x$-coordinates, and the Tate-Lichtenbaum and Ate pairings given both $x$-coordinates but at most one $y$-coordinate. This form of point compression is advantageous for elliptic curve pairings with small embedding degree, where one would be working over a field of large order (and consequently, taking a square root to recover $y$ could be expensive). The compression makes use of explicit recurrence formulas for elliptic curve point multiplication and for Miller functions in the case of embedding degree $k = 2$. As these recurrences are given solely in terms of the $x$-coordinate of the point, the pairings are also computed in terms of the $x$-coordinate of the points involved. Note, however, that neglecting the value of $y$ introduces a sign ambiguity, but this is resolved by taking the trace of the pairing, which is independent of the sign of $y$. It is perhaps worth investigating if the analogous results may be obtained for hyperelliptic pairings (for curves of the form $y^2 = F(x)$) of small embedding degree.

Another form of compression involves algebraic tori, which are $d$-dimensional generalizations of the multiplicative group $\G_m$. Naehrig, Barreto and Schwabe \cite{NBS:2008} use algebraic tori to compress computations, not just in the final exponentiation but also in the Miller loop of elliptic curve pairings. Their methods rely on explicit formulas for multiplication and squaring of torus elements and also exploit degree 6 twists. One might want to try similar methods for certain twists of hyperelliptic curves.

Another benefit of twists, as explained in Section \ref{ss:twistedate}, is that curves with a twist of degree $d$ allow one to use the twisted versions of Ate-type pairings. This means one computes the Miller function 
$f_{s,D_1}(D_2)$ for $D_1 \in \G_1$ and the divisor $D_2 = (u(x), v(x)) \in \G_2$ with $u(x)$ defined over the subfield $\F_{q^{k/(d,k)}}$, as opposed to computing $f_{s,D_2}(D_1)$.  Furthermore, the points of $\G_2$ can be represented as points on the Jacobian of the twist $\C'$ which allows for faster computations in the group $\G_2$. The example of Zhang \cite{Zhang:2008} uses a twist of degree $8$; to our knowledge,  pairings on curves with twists of degree 10 have not been implemented.

\subsection{Trace zero subvarieties}
\label{ss:tracezero}

For a hyperelliptic curve $\C$ of genus $g$ defined over $\F_q$, a
{\it trace zero subvariety} of $\C$ is a subgroup of the Jacobian of
$\C$ whose construction is connected to the Weil restriction of
scalars. The use of trace zero varieties for cryptographic
applications was first suggested by Frey \cite{Frey:2001}. The trace
zero subvariety of $\C$ over a field extension of degree $\ell$ is a
subgroup of $\Jac_C(\F_{q^\ell})$, which is isomorphic to the quotient
$\Jac_C(\F_{q^\ell})/ \Jac_C(\F_{q})$.

It can also be defined concretely as follows: Let $\pi$ be the $q^{\th}$
power Frobenius. Let $\ell$ be a prime and assume that $\ell
\nmid \#\Jac_C(\F_q)$.  We define the trace zero subvariety $G_\ell$ of
$\Jac_C(\F_{q^\ell})$ to be the set of elements of trace zero. I.e.,
\[
G_\ell(\F_q):= \{D \in \Jac_C(\F_{q^\ell}): D+ \pi(D)+ \dots+ \pi^{\ell-1}(D)= \mathcal{O}\}.
\]
Since $G_\ell(\F_q)$ is the kernel of the trace map, it is a subgroup of
$\Jac_C(\F_{q^\ell})$. To perform arithmetic in a trace zero subvariety
one can use the algorithms that work in the whole Jacobian. So far, no
specific algorithms for the group law are known that make use of the
subgroup properties.

Since $G_\ell(\F_q)$ is a subgroup of $\Jac_C(\F_{q^\ell})$, we can define a Tate-Lichtenbaum
pairing on it by restriction: suppose the order of
$G_\ell(\F_q)$ is divisible by a large prime factor $r$, but not by
$r^2$. Let $\G_1:= G_\ell[r] \cap \ker(\pi^\ell - [1])$ and $\G_2:=
G_\ell[r] \cap \ker(\pi^\ell-[q^\ell])$. Then the Tate-Lichtenbaum pairing on $G_\ell$ is a map
\[
t: \G_1 \times \G_2 \rightarrow \mu_{r}.
\]

On the points of $\G_1$, $\pi$ acts as multiplication by an integer
$s$ (\cite{DocheLangeCh:2006}),
and the same is true for the action of
$\pi$ on $\G_2$ (\cite[Proposition 3]{Cesena:2008}). Cesena 
\cite{Cesena:2008} gives a new algorithm for computing the Tate-Lichtenbaum
pairing over trace zero subvarieties of supersingular elliptic curves
by exploiting the action of the $q$-Frobenius. He uses the fact that
the $q$-Frobenius $\pi$ is an efficient endomorphism (rather than just
the $q^r$-Frobenius), together with
the fact that for particular supersingular elliptic curves the action of the
Frobenius can be computed more efficiently \cite[Lemmas 1--3]{Cesena:2008}. For these curves, the action of $\pi$ is (close to
being) multiplication of a power of $q$.

Experimentally, Cesena's algorithm is as efficient as the Tate-Lichtenbaum pairing on
supersingular elliptic curves, though less efficient than the eta
pairing $\eta_T$ or the optimal Ate pairing of Vercauteren.  It remains
to explore whether Cesena's algorithm generalizes to
supersingular hyperelliptic curves or non-supersingular
trace zero varieties. 

\subsection{Exploiting torsion groups of dimension $>2$}
\label{ss:highdimTorsion}
If $r$ is coprime with the characteristic of $\F_q$, the $r$-torsion
group of a Jacobian variety of dimension $g$ is isomorphic to
$(\Z/r\Z)^{2g}$. With the exception of the recent work by Okamoto and
Takashima \cite{OkamotoTakashima:2008}, all known pairing-based
cryptographic applications require only two linearly independent
torsion points and thus can be realized in the elliptic curve setting;
in fact, also the Okamoto-Takashima protocols can as well be implemented
using a product of two (supersingular) elliptic curves. It is an open
problem to find a cryptographic application that uses curves of genus
$2$ (or larger) and that does {\em not} work using elliptic curves. Both for
the ordinary and the supersingular case, constructions of Jacobians of
dimension $2$ with low {\em full} embedding degree (cf.\ Section
\ref{ss:ordinary}) are available (\cite{Freeman:2007,
  OkamotoTakashima:2008}).

\subsection{More Problems}
\label{ss:other}

For completeness, we include the  problems posed by Galbraith, Hess and Vercauteren  \cite{GHV:2007}, making note of any recent advancements:
\begin{enumerate}

\item {\it Construct pairing-friendly ordinary hyperelliptic curves with smaller $\rho$-values.} At this point in time, the smallest $\rho$-value obtained for an ordinary hyperelliptic curve of small embedding degree is $\rho \approx 20/9$ (for $g=2$, $k=27$; cf. Section \ref{s:friendly}). It is highly desirable to have curves with $\rho$-value $<2$. 

\item {\it Curves with $g \geq 3$.} For curves with $g \geq 3$, is it possible to develop efficient pairing-based cryptosystems which are also secure against the index calculus attacks available for these curves?

\item {\it Pairings on real models of hyperelliptic curves.} There have been recent examples \cite{GLM:2008} of efficient pairing computations on real models of hyperelliptic curves, as remarked in Section \ref{ss:looplength}. Are real models competitive with the imaginary models in general? Furthermore, are there efficient pairings on non-hyperelliptic curves?

\item {\it Torsion structure.} Is there an efficient method for selecting divisors  from $\Jac_{\C}(\F_{q^k})[r]$ for pairing computations? (See also Section \ref{ss:degeneratedivisors_faster}.) Furthermore, if this group has more than two generators, what cryptographic  applications are possible? (See also Section \ref{ss:highdimTorsion}.)

\item { \it Rubin-Silverberg point compression and Weil restriction.}  Can the Rubin-Silverberg method (see Section \ref{ss:rubinsilverberg}) be made 
more efficient in the elliptic curve case and/or generalized to Jacobians of curves of genus $g \geq 2$?

\item {\it Weil restriction.} As in Rubin-Silverberg point
  compression, certain abelian varieties can be identified with
  subvarieties of the Weil restriction of supersingular elliptic
  curves.  When the abelian variety is a Jacobian, are there
  explicitly computable homomorphisms between the elliptic curve and
  the Jacobian representation?
\end{enumerate}

%*******************************************************

\subsection*{Acknowledgments}
The authors are most grateful to David Freeman for his elaborate
feedback on earlier versions of this paper, which significantly
improved our work.  The authors further thank Paulo S.L.M.\ Barreto,
Anja Becker, Felix Fontein, Steven Galbraith, and Alfred Menezes for
helpful discussions and comments on an earlier draft of the paper.
Thanks to Alan Silvester for conscientious technical editing.  Thanks to
Jonathan Wise.  The work of the first author has
been supported by a National Science Foundation Graduate Research
Fellowship.  The third author has been supported by the Informatics Circle of Research Excellence Chair in Algorithmic Number Theory and Cryptography.
The fourth author is partially supported by National Science Foundation grant
DMS-0801123 and a Sloan Research Fellowship. The work of the fifth
author has been supported by a National Science Foundation Fellowship
(\#0802915) and a National Science and Engineering Research
Council Fellowship (\#373333). The fifth author is also grateful
to Harvard University and the Pacific Institute for the Mathematical
Sciences UBC, where some of this work was conducted.
Lastly, the sixth author acknowledges support by the National Science and Engineering
Research Council of Canada, and by the Premier's Research 
Excellence Award of the province of Ontario. 

The idea for this paper was conceived during the workshop ``Women in
Numbers'' at the Banff International Research Station in November 2008, where 
the authors formed a project group studying hyperelliptic pairings.

%********************************************************************

\providecommand{\bysame}{\leavevmode\hbox to3em{\hrulefill}\thinspace}
\providecommand{\MR}{\relax\ifhmode\unskip\space\fi MR }
\providecommand{\MRhref}[2]{%
  \href{http://www.ams.org/mathscinet-getitem?mr=#1}{#2}
}
\providecommand{\href}[2]{#2}


\begin{thebibliography}{10}

\bibitem{BalasubramanianKoblitz:1998}
R.~Balasubramanian and N.~Koblitz, \emph{The improbability that an elliptic
  curve has subexponential discrete log problem under the
  {M}enezes-{O}kamoto-{V}anstone algorithm}, Journal of Cryptology \textbf{11}
  (1998), no.~2, 141--145.

\bibitem{BGOS:2007}
P.S.L.M. Barreto, S.~Galbraith, C.~\'{O} h\'{E}igeartaigh, and M.~Scott,
  \emph{Efficient pairing computation on supersingular abelian varieties},
  Designs, Codes and Cryptography \textbf{42} (2007), 239--271.

\bibitem{BonehFranklin:2003}
D.~Boneh and M.~Franklin, \emph{Identity-based encryption from the {W}eil
  pairing}, SIAM Journal on Computing \textbf{32} (2003), no.~3, 586--615.

\bibitem{BLS:2004}
D.~Boneh, B.~Lynn, and H.~Shacham, \emph{Short signatures from the {W}eil
  pairing}, Journal of Cryptology \textbf{17} (2004), no.~4, 297--319.

\bibitem{BrezingWeng:2005}
F.~Brezing and A.~Weng, \emph{Elliptic curves suitable for pairing based
  cryptography}, Designs, Codes and Cryptography \textbf{37} (2005), 133--141.

\bibitem{Cantor:1987}
D.~Cantor, \emph{Computing in the {J}acobian of a hyperelliptic curve},
  Mathematics of Computation \textbf{48} (1987), no.~177, 95--101.

\bibitem{Cardona:2003}
G.~Cardona, \emph{On the number of curves of genus 2 over a finite field},
  Finite Fields and Their Applications \textbf{9} (2003), no.~4, 505--526.

\bibitem{Cardona:2005}
G.~Cardona, E.~Nart, and J.~Pujol{\`a}s, \emph{Curves of genus two over fields
  of even characteristic}, Mathematische Zeitschrift \textbf{250} (2005),
  no.~1, 177--201.

\bibitem{Cesena:2008}
E.~Cesena, \emph{Pairing with supersingular trace zero varieties revisited},
  Cryptology ePrint Archive Report 2008/404,
  \url{http://eprint.iacr.org/2008/404/}.

\bibitem{ChoieLee:2004}
Y.-J. Choie and E.~Lee, \emph{Implementation of {T}ate pairing on hyperelliptic
  curves of genus 2}, ICISC 2003, LNCS, vol. 2971, Springer-Verlag, 2004,
  pp.~97--111.

\bibitem{CocksPinch:2001}
C.~Cocks and R.G.E. Pinch, \emph{Identity-based cryptosystems based on the
  {W}eil pairing}, unpublished manuscript, 2001.

\bibitem{Coppersmith:1984}
D.~Coppersmith, \emph{Fast evaluation of logarithms in fields of characteristic
  two}, IEEE Transactions on Information Theory \textbf{30} (1984), 587--594.

\bibitem{DocheLangeCh:2006}
D.~Doche and T.~Lange, \emph{Handbook of {E}lliptic and {H}yperelliptic {C}urve
  {C}ryptography}, ch.~15 (Arithmetic of special curves), Chapman \& Hall/CRC,
  2006.

\bibitem{DuursmaLee:2003}
I.~Duursma and H.-S. Lee, \emph{Tate pairing implementation for hyperelliptic
  curves $y^2 = x^p-x+d$}, Advances in Cryptology - Asiacrypt 2003, LNCS, vol.
  2894, Springer-Verlag, 2003, pp.~111--123.

\bibitem{ELM:2004}
K.~Eisentr\"ager, K.~Lauter, and P.~Montgomery, \emph{Improved {W}eil and
  {T}ate pairings for elliptic and hyperelliptic curves}, Algorithmic Number
  Theory ANTS-VI, LNCS, vol. 3076, Springer-Verlag, 2004, pp.~169--183.

\bibitem{FGJ:2008}
X.~Fan, G.~Gong, and D.~Jao, \emph{Speeding up pairing computations on genus 2
  hyperelliptic curves with efficiently computable automorphisms}, Pairing
  2008, LNCS, vol. 5209, Springer-Verlag, 2008, pp.~243--264.

\bibitem{FreemanDataBase}
D.~Freeman, \emph{A generalized {B}rezing-{W}eng method for constructing
  pairing-friendly ordinary abelian varieties: Additional examples},
  \url{http://theory.stanford.edu/~dfreeman/papers/gen-bw-examples.pdf}.

\bibitem{Freeman:2007}
\bysame, \emph{Constructing pairing-friendly genus 2 curves with ordinary
  {J}acobians}, Pairing 2007, LNCS, vol. 4575, Springer-Verlag, 2007,
  pp.~152--176.

\bibitem{Freeman:2008}
\bysame, \emph{A generalized {B}rezing-{W}eng algorithm for constructing
  pairing-friendly ordinary abelian varieties}, Pairing 2008, LNCS, vol. 5209,
  Springer-Verlag, 2008, pp.~146--163.

\bibitem{FreemanSatoh}
D.~Freeman and T.~Satoh, \emph{Constructing pairing-friendly hyperelliptic
  curves using {W}eil restriction}, Preprint August 2009, 28 pages.

\bibitem{FST:2009}
D.~Freeman, M.~Scott, and E.~Teske, \emph{A taxonomy of pairing-friendly
  elliptic curves}, Journal of Cryptology, published online June 18, 2009. DOI:
  10.1007/s00145-009-9048-z (to appear in print).

\bibitem{FSS:2008}
D.~Freeman, P.~Stevenhagen, and M.~Streng, \emph{Abelian varieties with
  prescribed embedding degree}, Algorithmic Number Theory ANTS-VIII, LNCS, vol.
  5011, Springer-Verlag, 2008, pp.~60--73.

\bibitem{Frey:2001}
G.~Frey, \emph{Applications of arithmetical geometry to cryptographic
  constructions}, Finite Fields and Applications -- Proceedings of the Fifth
  International Conference on Finite Fields and Applications $F_q5$, Augsburg,
  August 2-6, 1999, Springer, 2001, pp.~128--161.

\bibitem{FreyLange:2006}
G.~Frey and T.~Lange, \emph{Fast bilinear maps from the {T}ate-{L}ichtenbaum
  pairing on hyperelliptic curves}, Algorithmic Number Theory ANTS-VII, LNCS,
  vol. 4076, Springer-Verlag, 2006, pp.~466--479.

\bibitem{FreyLangeCh23:2006}
\bysame, \emph{Handbook of {E}lliptic and {H}yperelliptic {C}urve
  {C}ryptography}, ch.~23 (Algebraic realizations of DL systems), Chapman \&
  Hall/CRC, 2006.

\bibitem{FreyLangeCh4:2006}
\bysame, \emph{Handbook of {E}lliptic and {H}yperelliptic {C}urve
  {C}ryptography}, ch.~4 (Background on curves and Jacobians), Chapman \&
  Hall/CRC, 2006.

\bibitem{FreyLangeCh5:2006}
\bysame, \emph{Handbook of {E}lliptic and {H}yperelliptic {C}urve
  {C}ryptography}, ch.~5 (Varieties over special fields), Chapman \& Hall/CRC,
  2006.

\bibitem{Galbraith:2001}
S.~Galbraith, \emph{Supersingular curves in cryptography}, Advances in
  Cryptology -- ASIACRYPT 2001, LNCS, vol. 2248, Springer-Verlag, 2001,
  pp.~495--513.

\bibitem{GHV:2007}
S.~Galbraith, F.~Hess, and F.~Vercauteren, \emph{Hyperelliptic pairings},
  Pairing 2007, LNCS, vol. 4575, Springer-Verlag, 2007, pp.~108--131.

\bibitem{GalbraithLin:2009}
S.~Galbraith and X.~Lin, \emph{Computing pairings using $x$-coordinates only},
  Designs, Codes and Cryptography \textbf{50} (2009), 305--324.

\bibitem{GLM:2008}
S.~Galbraith, X.~Lin, and D.~Mireles Morales, \emph{Pairings on hyperelliptic
  curves with a real model}, Pairing 2008, LNCS, vol. 5209, Springer-Verlag,
  2008, pp.~265--281.

\bibitem{GalbraithMcKeeValenca:2005}
S.~Galbraith, J.~McKee, and P.~Valen\c ca, \emph{Ordinary abelian varieties
  having small embedding degree}, Proc. Workshop on Mathematical Problems and
  Techniques in Cryptology, CRM, Barcelona, 2005, pp.~29--45.

\bibitem{GPRS:2009}
S.~Galbraith, J.~Pujol{\`a}s, C.~Ritzenthaler, and B.~Smith, \emph{Distortion
  maps for genus two curves}, Journal of Mathematical Cryptology \textbf{3}
  (2009), 1--18.

\bibitem{GalbraithScott:2008}
S.~Galbraith and M.~Scott, \emph{Exponentiation in pairing-friendly groups
  using homomorphisms}, Pairing 2008, LNCS, vol. 5209, Springer-Verlag, 2008,
  pp.~211--224.

\bibitem{GTTD:2007}
P.~Gaudry, E.~Thom{\'{e}}, N.~Th{\'{e}}riault, and C.~Diem, \emph{A double
  large prime variation for small genus hyperelliptic index calculus},
  Mathematics of Computation \textbf{76} (2007), 475--492.

\bibitem{GHOTV:2007}
R.~Granger, F.~Hess, R.~Oyono, N.~Th\'{e}riault, and F.~Vercauteren, \emph{Ate
  pairing on hyperelliptic curves}, Advances in Cryptology - Eurocrypt 2007,
  LNCS, vol. 4515, Springer-Verlag, 2007, pp.~419--436.

\bibitem{GPS:2006a}
R.~Granger, D.~Page, and N.~P. Smart, \emph{High security pairing-based
  cryptography revisited}, Algorithmic Number Theory ANTS-VII, LNCS, vol. 4076,
  Springer-Verlag, 2006, pp.~480--494.

\bibitem{GPS:2006b}
R.~Granger, D.~Page, and M.~Stam, \emph{On small characteristic algebraic tori
  in pairing-based cryptography}, LMS Journal of Computation and Mathematics
  \textbf{9} (2006), 64--85.

\bibitem{OhEigeartaighScott:2006}
C.~\'{O} h\'{E}igeartaigh and M.~Scott, \emph{Pairing calculation on
  supersingular genus 2 curves}, Selected Areas in Cryptography - SAC 2006,
  LNCS, vol. 4356, Springer-Verlag, 2007, pp.~302--316.

\bibitem{Hess:2008}
F.~Hess, \emph{Pairing lattices}, Pairing 2008, LNCS, vol. 5209,
  Springer-Verlag, 2008, pp.~18--38.

\bibitem{HSV:2006}
F.~Hess, N.~Smart, and F.~Vercauteren, \emph{The {E}ta pairing revisited}, IEEE
  Trans. Information Theory \textbf{52} (2006), 4595--4602.

\bibitem{Hitt:2007}
L.~Hitt, \emph{On the minimal embedding field}, Pairing 2007, LNCS, vol. 4575,
  Springer-Verlag, 2007, pp.~294--301.

\bibitem{Joux:2004}
A.~Joux, \emph{A one round protocol for tripartite {D}iffie-{H}ellman}, Journal
  of Cryptology \textbf{17} (2004), no.~4, 263--276.

\bibitem{KKAT:2004}
M.~Katagi, I.~Kitamura, T.~Akishita, and T.~Takagi, \emph{Novel efficient
  implementations of hyperelliptic curve cryptosystems using degenerate
  divisors}, WISA 2004, LNCS, vol. 3325, Springer-Verlag, 2005, pp.~345--359.

\bibitem{KawazoeTakahashi:2008}
M.~Kawazoe and T.~Takahashi, \emph{Pairing-friendly hyperelliptic curves of
  type $y^2 = x^5 + ax$}, Pairing 2008, LNCS, vol. 5209, Springer-Verlag, 2008,
  pp.~164--177.

\bibitem{Knuth:1997}
D.~Knuth, \emph{The {A}rt of {C}omputer {P}rogramming}, 3rd ed., vol.~2,
  Addison-Wesley, 1997.

\bibitem{Koblitz:89}
N.~Koblitz, \emph{Hyperelliptic cryptosystems}, Journal of Cryptology
  \textbf{1} (1989), no.~3, 139--150.

\bibitem{KoblitzMenezes:2005}
N.~Koblitz and A.~Menezes, \emph{Pairing-based cryptography at high security
  levels}, Cryptography and Coding 2005, LNCS, vol. 3796, Springer-Verlag,
  2005, pp.~13--36.

\bibitem{KohelDataBase}
D.~Kohel, \emph{Quartic {CM} fields database},
  \url{http://echidna.maths.usyd.edu.au/kohel/dbs/complex_multiplication2.html%
}.

\bibitem{Lange:2005}
T.~Lange, \emph{Formulae for arithmetic on genus 2 hyperelliptic curves},
  Applicable Algebra in Engineering, Communication and Computing \textbf{15}
  (2005), 295--328.

\bibitem{Lauter:2003}
K.~Lauter, \emph{The equivalence of the geometric and algebraic group laws for
  {J}acobians of genus 2 curves}, Topics in Algebraic and Noncommutative
  Geometry, Contemporary Mathematics, vol. 324, American Mathematical Society,
  2003, pp.~165--171.

\bibitem{LeeLeeLee:2008}
E.~Lee, H.-S. Lee, and Y.~Lee, \emph{Eta pairing computation on general
  divisors over hyperelliptic curves $y2=x^p-x+d$}, Journal of Symbolic
  Computation \textbf{43} (2008), 452--474.

\bibitem{LLP:2009}
E.~Lee, H.-S. Lee, and C.-M. Park, \emph{Efficient and generalized pairing
  computations on abelian varieties}, IEEE Transactions on Information Theory
  \textbf{55} (2009), 1793--1803.

\bibitem{LeeLee:2008}
E.~Lee and Y.~Lee, \emph{Tate pairing computation on the divisors of
  hyperelliptic curves of genus 2}, J. Korean Math. Soc. \textbf{45} (2008),
  no.~4, 1057--1073.

\bibitem{Lenstra:2001}
A.~K. Lenstra, \emph{Unbelievable security. {M}atching {AES} security using
  public key cryptosystems}, Advances in Cryptology - ASIACRYPT 2001, Lecture
  Notes in Computer Science, vol. 2248, Springer-Verlag, 2001, pp.~67--86.

\bibitem{LidlNiederreiter:1997}
R.~Lidl and H.~Niederreiter, \emph{Finite {F}ields}, 2nd ed., Cambridge
  University Press, 1997.

\bibitem{LMS:2004}
F.~Luca, D.~J. Mireles, and I.~E. Shparlinski, \emph{{MOV} attack in various
  subgroups on elliptic curves}, Illinois J. Math. \textbf{48} (2004),
  1041--1052.

\bibitem{Miller:2004}
V.~Miller, \emph{The {W}eil pairing and its efficient calculation}, Journal of
  Cryptology \textbf{17} (2004), 235--261.

\bibitem{Mumford:83}
D.~Mumford, \emph{Tata {L}ectures on {T}heta {I}, {II}}, Birkh{\"a}user,
  Boston, 1983/84.

\bibitem{NBS:2008}
M.~Naehrig, P.~Barreto, and P.~Schwabe, \emph{On compressible pairings and
  their computation}, AFRICACRYPT 2008, LNCS, vol. 5023, Springer-Verlag, 2008,
  pp.~371--388.

\bibitem{NIST:2006}
NIST, \emph{Recommendation for key management -- {P}art 1: General (revised)},
  NIST Special Publication 800-57, May 2006,
  \url{http://csrc.nist.gov/publications/nistpubs/800-57/SP800-57-Part1.pdf}.

\bibitem{Odlyzko:1985}
A.~Odlyzko, \emph{Discrete logarithms in finite fields and their cryptographic
  significance}, Advances in Cryptology - Eurocrypt 1984, LNCS, vol. 209,
  Springer-Verlag, 1985, pp.~224--314.

\bibitem{OkamotoTakashima:2008}
T.~Okamoto and K.~Takashima, \emph{Homomorphic encryption and signatures from
  vector decomposition}, Pairing 2008, LNCS, vol. 5209, Springer-Verlag, 2008,
  pp.~57--74.

\bibitem{PatersonCh:2005}
K.~Paterson, \emph{Advances in {E}lliptic {C}urve {C}ryptography}, ch.~X
  (Cryptography from pairings), Cambridge University Press, 2005.

\bibitem{Pollard:78}
J.~M. Pollard, \emph{Monte {C}arlo methods for index computation (mod $p$)},
  Mathematics of Computation \textbf{32} (1978), no.~143, 918--924.

\bibitem{RubinSilverberg:2002}
K.~Rubin and A.~Silverberg, \emph{Supersingular abelian varieties in
  cryptology}, Advances in Cryptology -- CRYPTO 2002, LNCS, vol. 2442,
  Springer-Verlag, 2002, pp.~336--353.

\bibitem{Satoh:2009}
T.~Satoh, \emph{The {B}rezing-{W}eng-{F}reeman method for certain genus 2
  hyperelliptic curves}, Advances in Cryptology -- Eurocrypt 2009, LNCS, vol.
  5479, Springer-Verlag, 2009, pp.~536--553.

\bibitem{Scott:2007}
M.~Scott, \emph{Implementing cryptographic pairings}, Pairing 2007, LNCS, vol.
  4575, Springer-Verlag, 2007, pp.~177--196.

\bibitem{Silverman:1992}
J.H. Silverman, \emph{The {A}rithmetic of {E}lliptic {C}urves}, Graduate Texts
  in Mathematics, vol. 106, Springer-Verlag, New York, 1992, Corrected reprint
  of the 1986 original.

\bibitem{OorschotWiener:99}
P.~C. van Oorschot and M.~J. Wiener, \emph{Parallel collision search with
  cryptanalytic applications}, Journal of Cryptology \textbf{12} (1999), 1--28.

\bibitem{Vercauteren:2008}
F.~Vercauteren, \emph{Optimal pairings}, Cryptology ePrint Archive Report
  2008/096, \url{http://eprint.iacr.org/2008/096/}.

\bibitem{Weng:2003}
A.~Weng, \emph{Constructing hyperelliptic curves of genus 2 suitable for
  cryptography}, Mathematics of Computation \textbf{72} (2003), 435--458.

\bibitem{WPP:2005}
T.~Wollinger, J.~Pelzl, and C.~Paar, \emph{Cantor versus {H}arley: Optimization
  and analysis of explicit formulae for hyperelliptic curve cryptosystem}, IEEE
  Transactions on Computers \textbf{54} (2005), 861--872.

\bibitem{Zhang:2008}
F.~Zhang, \emph{Twisted {A}te pairing on hyperelliptic curves and
  applications}, Cryptology ePrint Archive Report 2008/274,
  \url{http://eprint.iacr.org/2008/274/}.

\bibitem{ZZH:2008-2}
C.~Zhao, C.~Zhang, and J.~Huang, \emph{All pairings are in a group}, IEICE
  Transactions on Fundamentals of Electronics, Communications and Computer
  Sciences \textbf{E91-A} (2008), 3084--3087.

\bibitem{ZZH:2008-1}
\bysame, \emph{A note on the {A}te pairing}, International Journal of
  Information Security \textbf{7} (2008), 379--382.

\end{thebibliography}
\end{document}